\documentclass[reqno]{amsart}
\usepackage{amsmath,amsthm,xspace,enumerate,amscd,graphicx,epstopdf}
\usepackage{mathrsfs,float}
\usepackage{amssymb,url,color}
\usepackage{subfig}
\usepackage[pdftex]{hyperref}
\usepackage{array,multirow,makecell,booktabs,diagbox}
\usepackage[format=hang,font=small]{caption}

\newtheorem{theorem}{Theorem}[section]
\newtheorem{remark}{Remark}[section]
\newtheorem{lemma}{Lemma}[section]
\newtheorem{proposition}{Proposition}[section]

\begin{document}

\author{Nannan Ma}
\address[N. Ma]{Zhengzhou Zhongyuan Sub-branch, Agricultural Bank of China, 450000, Henan,
China.}
\email{\href{mailto: N. Ma <1055165454@qq.com>}{1055165454@qq.com}}

\author{Hailin Sang}
\address[H. Sang]{ Department of Mathematics, University of Mississippi, University,  MS, 38677,
USA.}
\email{\href{mailto: H. Sang<sang@olemiss.edu
>}{sang@olemiss.edu}}

\author{Guangyu Yang}
\address[G. Yang]{ School of Mathematics and Statistics, Zhengzhou University, 450001, Henan,
China.}
\email{\href{mailto: G. Yang
<guangyu@zzu.edu.cn>}{guangyu@zzu.edu.cn}}




\begin{abstract}
We establish the asymptotic theory of least absolute deviation estimators for AR(1) processes with autoregressive parameter satisfying $n(\rho_n-1)\to\gamma$ for some fixed $\gamma$ as $n\to\infty$,
which is parallel to the results of ordinary least squares estimators developed by Andrews and Guggenberger (2008) in the case $\gamma=0$ or Chan and Wei (1987) and Phillips (1987) in the case $\gamma\ne 0$.  Simulation experiments are conducted to confirm the theoretical results and to demonstrate the robustness of the least absolute deviation estimation.
\end{abstract}

\title[LAD estimation for AR(1) processes]
{Least absolute deviation estimation for AR(1) processes with roots close to unity}

\maketitle


\section[Introduction]{Introduction}

\noindent  Consider the following AR(1) process
\begin{align*}
y_i=\rho y_{i-1}+\epsilon_i, \quad 1\leq i\leq n,
\end{align*}
where $\rho$ is a deterministic parameter and $\{\epsilon_i\}_{i\in\mathbb{Z}}$ is a sequence of independent and identically distributed (i.i.d.) random variables with mean zero and finite variance. The asymptotic properties of the ordinary least squares (OLS) estimator of $\rho$
have been extensively studied in the literature; we refer to Anderson (1959), White (1958), Dickey and Fuller (1979), Phillips (1987), Chan (2009), Miao and Shen (2009), and references therein.
To further handle the data that allows for large shocks in dynamic structure of the process, e.g. modeling the asset-price bubbles, it is usually to require the parameter $\rho$ to depend on the sample size $n$. On the other hand, to understand the phenomena that the unit root test  generally has a low discriminatory power against the alternative of root close to but not equal to unity, Bobkoski (1983) and Cavanagh (1985) introduced a local unit root model with the parameter $\rho$ depending on the sample size $n$ and tending to unity as $n\to\infty$.    Since then, many researchers systematically established the limit theory for various near unit root processes. See Chan and Wei (1987), Phillips (1987, 1988), Phillips and Magdalinos (2007), Aue and Horv\'{a}th (2007), Andrews and Guggenberger (2008), Buchmann and Chan (2013), Miao {\it et al.} (2015), Jiang {\it et al.} (2022), Tanaka (2017), and references therein. In particular, we refer to Stock (1991) for the empirical research or Phillips (2021) for the recent theoretical progress and empirical research on the processes with near unit roots.

\vskip5pt
Recently, Zhou and Lin (2014) and Wang {\it et al.} (2020) studied the statistical inference for the autoregressive parameter under the framework of the least absolute deviation (LAD) estimation. They proved that, if $\rho_n\to1$ and $n|\rho_n-1|\to\infty$ as $n\to\infty$, then the LAD estimators have normal and Cauchy asymptotic distributions under the mildly-stationary case and the mildly-explosive case respectively, which are complementary to the results of the OLS estimators established by Giraitis and Phillips (2006) and Phillips and Magdalinos (2007). In fact, the LAD estimator was first considered in the study of  autoregressive processes in the case that the regressive parameters are constants and the innovations have infinite variance due to the robust property of LAD estimation. For example, see the papers, Pollard (1991), Phillips (1991), Davis {\it et al.} (1992), and Li and Li (2009). Specially, Herce (1996) studied the asymptotic property of the LAD estimator in a unit root process with finite variance innovations and correspondingly developed the unit root test in this case which complement similar ones obtained by Knight (1989).
\vskip5pt

Motivated by the above work, the goal of the present article is to establish the asymptotic theory of LAD estimators for AR(1) processes when the regressive parameter satisfies $n(\rho_n-1)\to0$ or $n(\rho_n-1)\to\gamma$ for some fixed $\gamma\ne 0$ as $n\to\infty$. To the best of our knowledge, this part of research is still missing although we already have a rich literature in LAD estimations and unit root test for AR(1) processes. It is shown that, if $n(\rho_n-1)\to0$ as $n\to\infty$,  the limiting distributions are dominated by the initial conditions both in the near-stationary case and the near-explosive case. This phenomena have been studied by M\"{u}ller and Elliott (2003), Phillips and Magdalinos (2009), and references therein. Our results in the near-stationary case correspond exactly to the theory on the OLS estimators for the same model developed by Andrews and Guggenberger (2008). With the condition $n(\rho_n-1)\to\gamma$ for some fixed $\gamma\ne 0$ as $n\to\infty$, we study the asymptotic theory of the LAD estimator under the assumption $y_0=0$. The work in this case is complementary to the theory of the OLS estimation developed by Chan and Wei (1987) and Phillips (1987).
\vskip5pt

This paper is organized as follows. After introducing the model, we state the main results in Section \ref{s2}. Section \ref{s4} reports some simulation studies to illustrate the finite sample performance and to confirm the asymptotic results. Section \ref{s5} consists of concluding remarks and Section \ref{s6} provides the proofs of main results. All the other technical proofs are given in Appendix.
\vskip5pt

Throughout this paper, the symbols `$\rightarrow_d$', `$\rightarrow_p$' and `$=_d$' denote the weak convergence, convergence in probability and equality in distribution, respectively; $o_p(1)$ means tending to zero in probability; $\mathcal{C}$ denotes the standard Cauchy random variable and $\mathscr{N}(\mu,\sigma^2)$ denotes the normal random variable with  mean $\mu$ and variance $\sigma^2$; $[x]$ represents the integral part of $x$; ${\rm sign}(\cdot)$ is the signum function; $\mathbb{I}_{A}$ is the indicator function of set $A$ and $\det(M)$ is the determinant of matrix $M$.



\section[main results]{Main results}\label{s2}

Suppose that the data $\{y_i\}$ are generated by the following autoregressive model with order one
\begin{align}\label{model}
y_i=\rho_n y_{i-1}+\epsilon_i, \quad 1\leq i\leq n,
\end{align}
where the parameter $\rho_n$ satisfies $\rho_n\to 1$ as $n\to\infty$, and the noises $\{\epsilon_i\}_{i\in\mathbb{Z}}$ satisfy the following assumptions:
\begin{itemize}
\item[(i).] $\{\epsilon_i\}_{i\in\mathbb{Z}}$ are i.i.d. random variables with mean zero and finite variance $\sigma^2$;
\item[(ii).] $\epsilon_1$ has zero median and a differentiable density function $f(x)$ in $\mathbb{R}$ with $f(0)>0$ and $\sup_{x\in\mathbb{R}}|f'(x)|<\infty$.
\end{itemize}
Note that, formally the data $\{y_{i}\}$ is a triangular array but here $n$ is omitted for notational simplicity. With the above assumptions on the noises, one can also study the estimation of the parameter $\rho_n$ with $\rho_n\to -1$ as $n\to\infty$. For the purpose of simplicity, we only state the results for the positive $\rho_n$ case in this paper.
\vskip5pt

The LAD estimator of $\rho_n$ is defined as a solution of the following extremum problem
\begin{align}\label{LAD est}
\hat{\rho}_{\rm LAD}=\arg\min_{\rho\in\mathbb{R}}\left\{\frac{1}{n}\sum_{i=1}^n|y_i-\rho y_{i-1}|\right\}.
\end{align}
Notice that, the LAD estimators usually do not have closed forms and are not unique if the object function has a flat segment; see, e.g. Herce (1996). Moreover, the LAD estimators are robust, especially they are not significantly affected by the presence of outliers. This is confirmed by the simulation study in Section \ref{s4}.


\subsection[The case $n(\rho_n-1)\to0$ as $n\to\infty$]{The case $n(\rho_n-1)\to0$ as $n\to\infty$}

We first consider the case that the parameter $\rho_n$ satisfies $n(\rho_n-1)\to0$ as $n\to\infty$. That is, the autoregressive parameter $\rho_n$ is very nearly unity in the sense that $\rho_n$ is away from unity by $o(n^{-1})$. The following two theorems are the  main results in this case. One is for the near-stationary case $0<\rho_n<1$ and the other is for the near-explosive case $\rho_n>1$.

\begin{theorem}[The near-stationary case]\label{near-s0}
For model (\ref{model}) with $0<\rho_n<1$ and $n(\rho_n-1)\to0$ as $n\to\infty$, assume that $y_0$ depends on the full past of the noise, i.e. $y_0=\sum_{j=0}^{\infty}\rho_n^j\epsilon_{-j}$, then, as $n\to\infty$, we have
\begin{align}
\sqrt{\frac{n}{1-\rho_n^2}}\,(\hat{\rho}_{\rm LAD}-\rho_n)\rightarrow_d\frac{\mathcal{C}}{2\sigma f(0)}
\end{align}
and
\begin{align}
\sqrt{\sum_{i=1}^ny_{i-1}^2}\,(\hat{\rho}_{\rm LAD}-\rho_n)\rightarrow_d\mathscr{N}\left(0,\frac{1}{4f^2(0)}\right).
\end{align}
\end{theorem}

\begin{theorem}[The near-explosive case]\label{near-e0}
For model (\ref{model}) with $\rho_n>1$ and $n(\rho_n-1)\to0$ as $n\to\infty$, assume that $y_0$ is an infinitely distant initialization, i.e. $y_0=\sum_{j=0}^{\kappa_n}\rho_n^j\epsilon_{-j}$, with $\kappa_n/n\to\infty$ and $\kappa_n(\rho_n-1)\to0$ as $n\to\infty$, then, as $n\to\infty$, we have
\begin{align}
\sqrt{n\kappa_n}\,(\hat{\rho}_{\rm LAD}-\rho_n)\rightarrow_d\frac{\mathcal{C}}{2\sigma f(0)}
\end{align}
and
\begin{align}
\sqrt{\sum_{i=1}^ny_{i-1}^2}\,(\hat{\rho}_{\rm LAD}-\rho_n)\rightarrow_d\mathscr{N}\left(0,\frac{1}{4f^2(0)}\right).
\end{align}
\end{theorem}

We give some comments on Theorems \ref{near-s0} and \ref{near-e0}.

\begin{remark}\label{verynear1}
Since the AR(1) model (\ref{model}) is causal  when $0<\rho_n<1$, the initial value $y_0$ can be written as a linear combination over the past information in the linear process form $y_0=\sum_{j=0}^{\infty}\rho_n^j\epsilon_{-j}$. In the same case, $\rho_n\in(-1,1),\;n(\rho_n-1)\to0$ and $y_0$ depends on the full past of the noise, Andrews and Guggenberger (2008) obtained the limit theorems for the OLS estimators which are similar to Theorem \ref{near-s0} except for the appearance of $f(0)$ here. Furthermore, because the limiting distributions of the $t$-type estimator include $f(0)$ in LAD estimations, in practice we use the following density estimator (Silverman, 1986) for statistical inference
\[
\hat{f}_n(0):=\frac{1}{nb_n}\sum_{i=1}^nK\left(\frac{y_i-\hat\rho_{\rm LAD}y_{i-1}}{b_n}\right),
\]
where $b_n$ is the bandwidth and $K(\cdot)$ is a kernel function, e.g. the Gaussian or logistic kernel.
\end{remark}

\begin{remark}\label{verynear2}
As we mentioned previously, Zhou and Lin (2014) and Wang {\it et al.} (2020) studied the LAD estimations for $\rho_n$ satisfying $n|\rho_n-1|\to\infty$ where the initial value $y_0=o_p(|\rho_n-1|^{-1/2})$ independent of $\{\epsilon_i, i\geq1\}$. To be explicit, they obtained that
\begin{align}\label{zhoulinwang1}
\left\{
\begin{array}{ll}
\sqrt{\frac{n}{1-\rho_n^2}}\,(\hat{\rho}_{\rm LAD}-\rho_n)\rightarrow_d\mathscr{N}\left(0,\frac{1}{4f^2(0)}\right)\\
\sqrt{\sum_{i=1}^ny_{i-1}^2}\,(\hat{\rho}_{\rm LAD}-\rho_n)\rightarrow_d\mathscr{N}\left(0,\frac{1}{4f^2(0)}\right)
\end{array}
\right.
\end{align}
and
\begin{align}\label{zhoulinwang2}
\frac{\rho_n^n}{\rho_n^2-1}\,(\hat{\rho}_{\rm LAD}-\rho_n)\rightarrow_d\frac{\mathcal{C}}{2\sigma f(0)},
\end{align}
for the mildly-stationary case and the mildly-explosive case, respectively. Their results match the classic results in Phillips and Magdalinos (2007) except for the appearance of $f(0)$, just like the usual LAD estimations.
\vskip5pt

It is worth pointing out that, here for the very nearly unit root processes, Theorem \ref{near-s0} and Theorem \ref{near-e0} show that, in both near-stationary and near-explosive cases, the initial value $y_0$ dominates the asymptotic distributions of the LAD estimator and the $t$-type estimator which are Cauchy and normal respectively. Moreover, since $\rho_n\rightarrow 1$ at a much faster rate than $1/n$, for the near-stationary case, by Lemma \ref{l-s0}, the assumption $y_0=o_p(|\rho_n-1|^{-1/2})$ can not hold and the convergence rate $\sqrt{n/(1-\rho_n^2)}$ has a larger order than $n$ which enlarges the convergence rate spectrum, while, for the near-explosive case, the convergence rate $\sqrt{n\kappa_n}$ also has a large order than $n$ which is different from that in equation (\ref{zhoulinwang2}).
\end{remark}



\subsection{The case $n(\rho_n-1)\to\gamma$ as $n\to\infty$}

Now we consider the local unit root case that the parameter $\rho_n$ satisfies $n(\rho_n-1)\to\gamma$ for some fixed $\gamma\ne 0$ as $n\to\infty$.
We remark that, unlike the case, $\rho_n\ne1$ and $n(\rho_n-1)\to0$ as $n\to\infty$, where the initial condition entirely dominates the asymptotic distribution of the sample variance $\sum_{i=1}^ny_{i-1}^2$, in this case, the initial value depending on the past of the noise does not affect the analysis methods essentially so we can assume that $y_0=0$ for simplicity.

\vskip5pt

To state the main results, we first introduce some notations. For $1\leq i\leq n$, let
\[
K_{n,i}:=\frac{1}{\sqrt{n}}\rho_n^{i-1-n}{\rm sign}(\epsilon_i)\quad{\text{and}}\quad L_{n,i}:=\frac{1}{\sqrt{n}}\rho_n^{n-i}\epsilon_i,
\]
and define
\begin{align}\label{kn}
K_n(t):=\sum_{i=1}^{[nt]}K_{n,i} \quad{\text{and}}\quad L_n(t):=\sum_{i=1}^{[nt]}L_{n,i},
\end{align}
for $0\leq t\leq1$. If we further assume that $\sigma^2=1$ and $\mathbb{E}|\epsilon_1|^{2+\delta}<\infty$ for some $\delta>0$, by Lemma \ref{KLM1}, the process $\{(K_n(t),L_n(t)), 0\leq t\leq1\}$ converges weakly to a continuous process
\begin{align}\label{xt}
{\bf{X}}(t)=(K(t),L(t)), \quad 0\leq t\leq1
\end{align}
with independent Gaussian increments, mean vector zero and covariance matrix
\begin{align}\label{covm}
\Gamma(t)=(\gamma_{ij}(t)):=\left(
                                  \begin{array}{cc}
                                    \frac{e^{-2\gamma}}{2\gamma}(e^{2\gamma t}-1) & t\mathbb{E}|\epsilon_1| \\
                                    t\mathbb{E}|\epsilon_1| & \frac{e^{2\gamma}}{2\gamma}(1-e^{-2\gamma t}) \\
                                  \end{array}
                                \right), \quad 0\leq t\leq1.
\end{align}

\begin{theorem}\label{main result2}
 For model (\ref{model}), assume that  $y_0=0$, $\sigma^2=1$, $n(\rho_n-1)\to\gamma$ for some fixed $\gamma\neq0$ as $n\to\infty$, and $\mathbb{E}|\epsilon_1|^{2+\delta}<\infty$ for some $\delta>0$, then, as $n\to\infty$, we have
\begin{align}
n(\hat{\rho}_{\rm LAD}-\rho_n)&\rightarrow_d\frac{1}{2f(0)}\cdot\frac{\int_0^1L(t)\,{\rm d}K(t)}{\int_0^1e^{-2\gamma(1-t)}L^2(t)\,{\rm d}t}
\end{align}
and
\begin{align}
\sqrt{\sum_{i=1}^ny_{i-1}^2}\,(\hat{\rho}_{\rm LAD}-\rho_n)&\rightarrow_d\frac{1}{2f(0)}\cdot\frac{\int_0^1L(t)\,{\rm d}K(t)}{\sqrt{\int_0^1e^{-2\gamma(1-t)}L^2(t)\,{\rm d}t}}.
\end{align}
\end{theorem}


\begin{remark}
Zhou and Lin (2014) and Wang {\it et al.} (2020) worked on the case that $\rho_n\to1$ and $n|\rho_n-1|\to\infty$ as $n\to\infty$, where the convergence rate for the  asymptotic distributions of the LAD estimators has an order in the range $(n^{1/2}, n)$ for the mildly-stationary case. Together with the results in Theorem \ref{near-s0} and Theorem \ref{near-e0}, we enlarge the convergence rate spectrum in the asymptotic distributions of the LAD estimators. The similar relationship between the rate of $\rho_n$ and the rate in the asymptotic distributions was observed for OLS estimators in the literature, see, e.g. Chan and Wei (1987), Phillips (1987), Phillips and Magdalinos (2007), and Andrews and Guggenberger (2008).
\end{remark}

\begin{remark}\label{DL}
Let
\begin{align}
\mathscr{D}(\gamma):=\frac{1}{2f(0)}\cdot\frac{\int_0^1L(t)\,{\rm d}K(t)}{\int_0^1e^{-2\gamma(1-t)}L^2(t)\,{\rm d}t}
\end{align}
and
\begin{align}
\mathscr{L}(\gamma):=\frac{1}{2f(0)}\cdot\frac{\int_0^1L(t)\,{\rm d}K(t)}{\sqrt{\int_0^1e^{-2\gamma(1-t)}L^2(t)\,{\rm d}t}}.
\end{align}
Notice that $\mathscr{D}(\gamma)$ and $\mathscr{L}(\gamma)$ are continuous families of distributions indexed by the parameter $\gamma$.  Thus, for the case, $\rho_n=1+{\gamma}/{n}$, by Lemma \ref{KLM1}, when $\gamma=0$, i.e. $\rho_n\equiv1$, we have
\begin{align*}
\mathscr{D}(0)=\frac{1}{2f(0)}\cdot\frac{\int_0^1W_2(t)\,{\rm d}W_1(t)}{\int_0^1 W_2^2(t)\,{\rm d}t}\quad{and}\quad \mathscr{L}(0)=\frac{1}{2f(0)}\cdot\frac{\int_0^1W_2(t)\,{\rm d}W_1(t)}{\sqrt{\int_0^1 W_2^2(t)\,{\rm d}t}},
\end{align*}
here ${\bf W}=(W_1,W_2)$ is a bivariate Brownian motion with covariance matrix
\begin{align*}
\Xi=\left(
  \begin{array}{cc}
    1 & \mathbb{E}|\epsilon_1| \\
    \mathbb{E}|\epsilon_1| & 1\\
  \end{array}
\right).
\end{align*}
 These exactly correspond to the main results in Herce (1996), where the author discussed the structure of $\mathscr{D}(0)$, a combination of a ``unit root'' distribution and a scale mixture of normal distributions, and used it to construct the LAD-based unit root tests. In addition, since $\sigma^2=1$, we have $\mathbb{E}|\epsilon_1|\leq1$, hence the matrix $\Xi$ is non-negative definite. However, here $\mathscr{D}(\gamma)$ and $\mathscr{L}(\gamma)$ are too complicate to be analyzed effectively so we do some simulations to give an overall view on the distributions of $\mathscr{D}(\gamma)$ and $\mathscr{L}(\gamma)$ in Section \ref{s4} which illustrate how they depend on the parameter $\gamma$. Finally, we also remark that, from Theorem 1 in Herce (1996), the local unit root case has the optimal rate of convergence for the alternative hypothesis, that is, it has the same rate as in the unit root case.
\end{remark}

Followed the above remark, there is a natural question, what can we say about the distributions of $\mathscr{L}(\gamma)$ as $|\gamma|\to\infty$?
In fact, like the case for the OLS estimators in Chan and Wei (1987) and  Phillips (1987), we have the following result.

\begin{theorem}\label{main result3}
For model (\ref{model}), assume that $y_0=0$, $\sigma^2=1$, $n(\rho_n-1)\to\gamma$ for some fixed $\gamma\neq0$ as $n\to\infty$, and $\mathbb{E}|\epsilon_1|^{2+\delta}<\infty$ for some $\delta>0$, then we have
\begin{align}
2f(0)\cdot\mathscr{L}(\gamma)\rightarrow_d \mathscr{N}(0,1), \quad as \; |\gamma|\to\infty.
\end{align}
\end{theorem}

\begin{remark}
Recall from Remark \ref{verynear2}, Zhou and Lin (2014) and Wang {\it et al.} (2020) proved that
\begin{align}\label{zhoulinwang3}
2f(0)\cdot\sqrt{\sum_{i=1}^ny_{i-1}^2}\,(\hat{\rho}_{\rm LAD}-\rho_n)\rightarrow_d\mathscr{N}\left(0,1\right),
\end{align}
for the mildly-stationary case, $0<\rho_n<1$ and $n(1-\rho_n)\to\infty$ as $n\to\infty$. Theorem \ref{main result3} together with Theorem \ref{main result2} yields that, if let $n(\rho_n-1)\to\gamma$ first, and then let $|\gamma|\to\infty$, we can also get the asymptotic normal distribution.  Moreover, it is worth noting that, although Theorem \ref{main result3} holds for $\gamma\to\infty$ and $\gamma\to-\infty$, the underlying reasoning is quite different between these two cases. In fact, for the large negative $\gamma$'s, $\rho_n$ can be thought of much less than one so the model (\ref{model}) can
be regarded as stationary and then the asymptotic normality holds; for the positive $\gamma$, $\rho_n$ is larger than one and the model (\ref{model}) is explosive, however Theorem \ref{main result3} shows that the asymptotic normality is still valid which is different from the mildly-explosive case discussed in Zhou and Lin (2014) and Wang {\it et al.} (2020), i.e. the equation (\ref{zhoulinwang2}) in Remark \ref{verynear2}.
\end{remark}





\section{simulations}\label{s4}

In this section, we work on Monte Carlo simulation to examine the finite sample performance of the estimators in our main results. For convenience, in all experiments, we always suppose that the data are generated by the AR(1) model, $y_i=\rho_n y_{i-1}+\epsilon_i$ with $\rho_n=1+\gamma n^{-\beta}$ for some $\gamma\in\mathbb{R}$ and $\beta\geq1$. Furthermore, we also assume that the innovations are i.i.d. and  have $\mathscr{N}(0,1)$ or $U(-1,1)$ distributions. Here $U(-1,1)$ denotes the uniform random variable on the interval $(-1, 1)$.
To estimate the parameter $f(0)$, we apply the following density estimator (Silverman, 1986)
\[
\hat f_n(0)=\frac{1}{nb_n}\sum_{i=1}^nK\left(\frac{y_i-\hat\rho_{\rm LAD}y_{i-1}}{b_n}\right),
\]
where $K(\cdot)$ is the Gaussian kernel and the optimal bandwidth $b_n$ associated with the Matlab function {\it ksdensity}($\cdot,0$) is automatically selected from the data.
\vskip5pt

\noindent\textbf{Density curves of estimates.} We first consider the simulations of Theorem \ref{near-s0} and Theorem \ref{near-e0}. For the near-stationary case, $y_0=\sum_{j=0}^{\infty}\rho_n^{j}\epsilon_{-j}$, the true parameters are taken to be $(\gamma,\beta)=(-10,2)$  and we denote the normalized estimator by $\hat{\rho}_{s}=2\sigma\hat f_n(0)\sqrt{n/(1-\rho_n^2)}(\hat{\rho}_{\rm LAD}-\rho_n)$. For the near-explosive case, take $y_0=\sum_{j=0}^{\kappa_n}\rho_n^{j}\epsilon_{-j}$, $(\gamma,\beta)=(10,2)$, $\kappa_n=[n^{1.3}]$, and denote the normalized estimator by $\hat{\rho}_{e}=2\sigma \hat f_n(0)\sqrt{n\kappa_n}(\hat{\rho}_{\rm LAD}-\rho_n)$. For both cases, denote the $t$-type estimators by $T_n:=2\hat{f}_n(0)\sqrt{\sum_{i=1}^ny_{i-1}^2}\,(\hat{\rho}_{\rm LAD}-\rho_n)$. We simulate $3000$ replications with sample size $n=1000$ for each case. Figure 1 shows that the density curves of $\hat{\rho}_{s}$ and $\hat{\rho}_{e}$ are close to that of the standard Cauchy random variable.
Figure 2 confirms the asymptotic normality of $T_n$ in the near-stationary and near-explosive cases by using the Q-Q graphs.
\vskip5pt

Next we simulate the limiting distributions, $\mathscr{D}(\gamma)$ and $\mathscr{L}(\gamma)$, defined as in Remark \ref{DL}.
Let $\beta=1$ and suppose that the innovations are standard normal random variables so $f(0)=1/\sqrt{2\pi}$. For $\mathscr{D}(\gamma)$, simulate $1000$ replications with sample size $n=200$ for each case; Figure 3 illustrates the shape of the asymptotic density curves of $2f(0)\mathscr{D}(\gamma)$ for different $\gamma$. For $\mathscr{L}(\gamma)$, simulate $1000$ replications with sample size $n=2000$ ($\gamma<0$) and $n=1000$ ($\gamma>0$) for each case; Figure 4 confirms Theorem \ref{main result3}, i.e. the asymptotic normality of $2f(0)\mathscr{L}(\gamma)$ as $|\gamma|\to\infty$, which also demonstrates the difference of the convergence rate between $\gamma\to\infty$ and $\gamma\to-\infty$.
\begin{figure}[H]
\subfloat[The near-stationary case]{
\begin{minipage}[t]{0.5\linewidth}
\centering
\includegraphics[width=6.5cm,height=6.5cm]{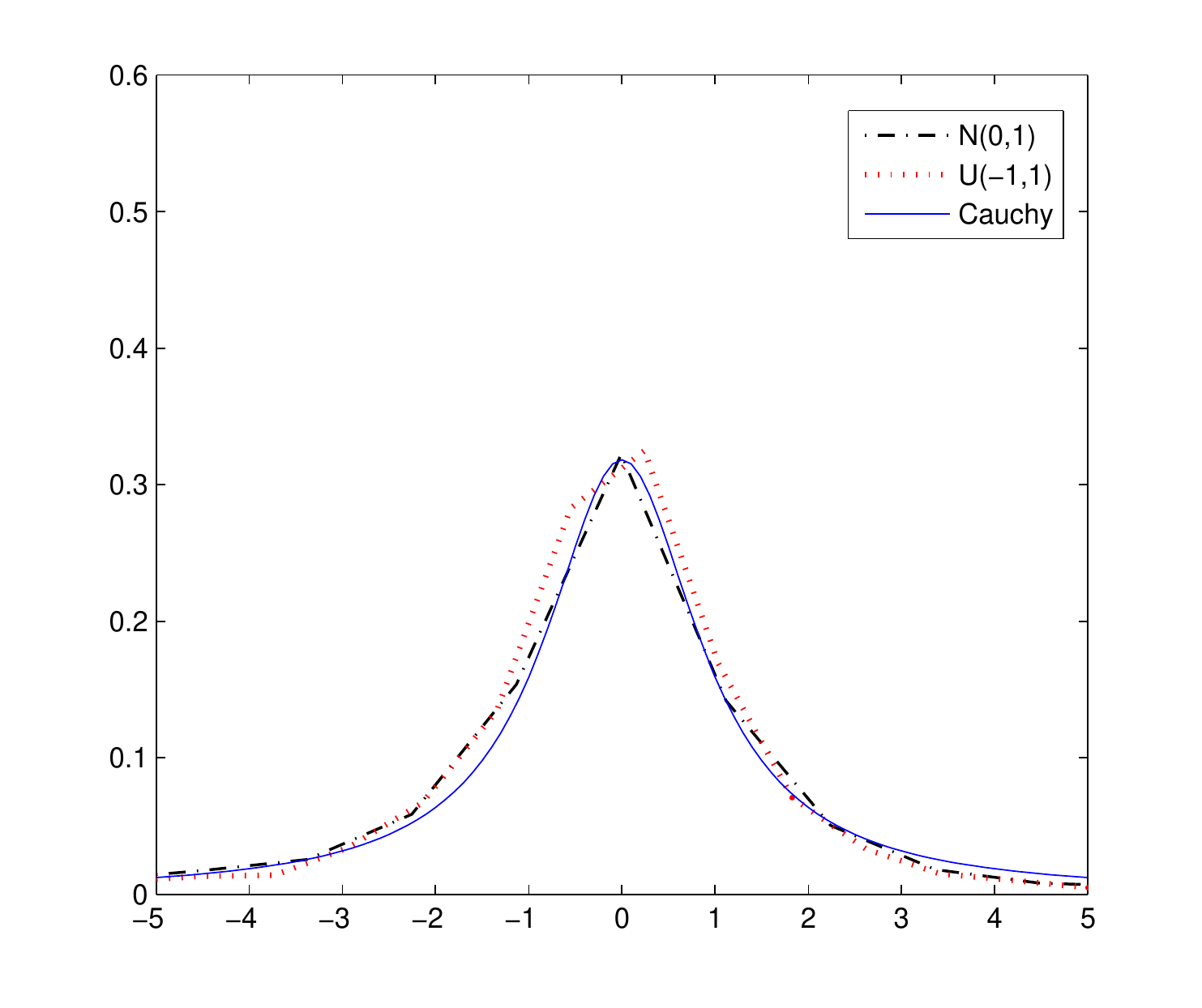}
\end{minipage}
}
\subfloat[The near-explosive case]{
\begin{minipage}[t]{0.5\linewidth}        
\includegraphics[width=6.5cm,height=6.5cm]{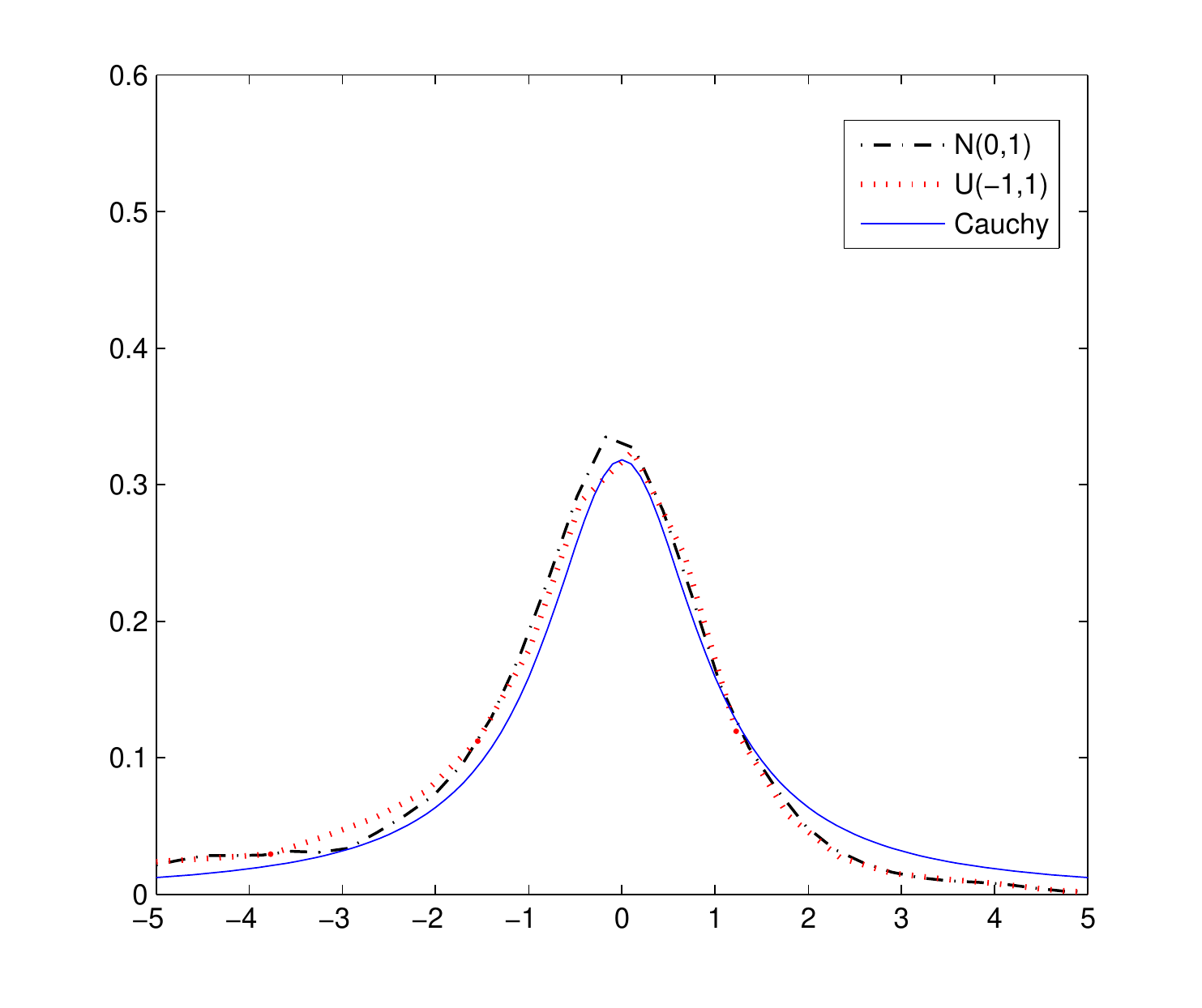}
\end{minipage}
}
\caption{\small{Density curves of $\hat{\rho}_{s}$ and $\hat{\rho}_{e}$ with $\mathscr{N}(0,1)$ and $U(-1,1)$ innovations.}}
\end{figure}

\begin{figure}[H]
\subfloat[The near-stationary case]{
\begin{minipage}[t]{0.5\linewidth}
\centering
\includegraphics[width=6.5cm,height=6.5cm]{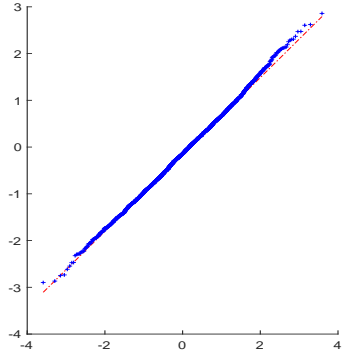}
\end{minipage}
}
\subfloat[The near-explosive case]{
\begin{minipage}[t]{0.5\linewidth}        
\includegraphics[width=6.5cm,height=6.5cm]{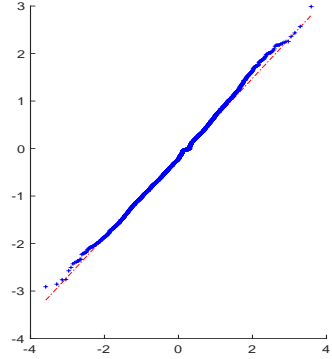}
\end{minipage}
}
\caption{\small{Q-Q graphs of statistics $T_n$ with $\mathscr{N}(0,1)$ innovation.}}
\end{figure}


\begin{figure}[H]
\subfloat[The near-stationary case]{
\begin{minipage}[t]{0.5\linewidth}
\centering
\includegraphics[width=6.5cm,height=6.5cm]{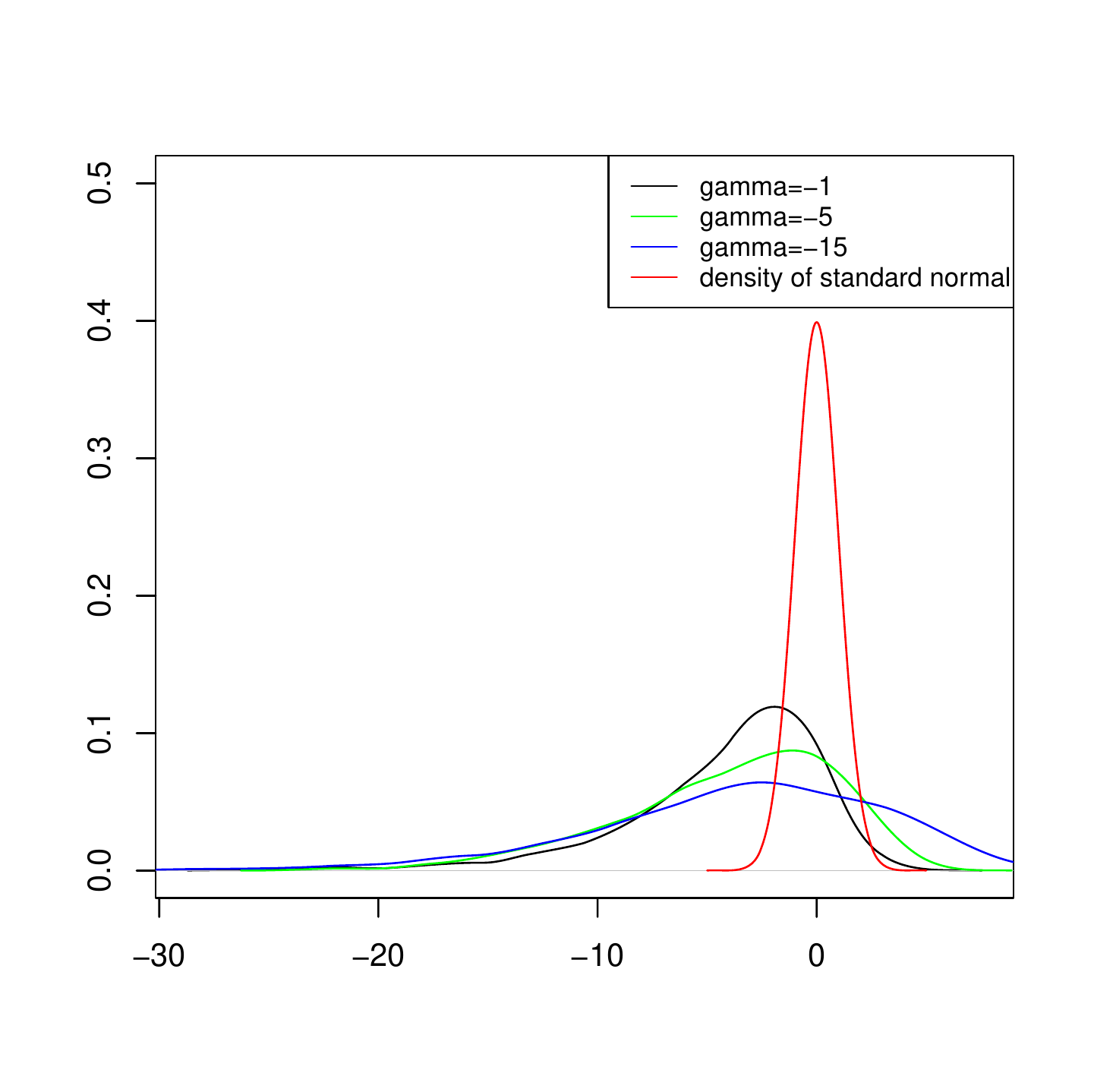}
\end{minipage}
}
\subfloat[The near-explosive case]{
\begin{minipage}[t]{0.5\linewidth}        
\includegraphics[width=6.5cm,height=6.5cm]{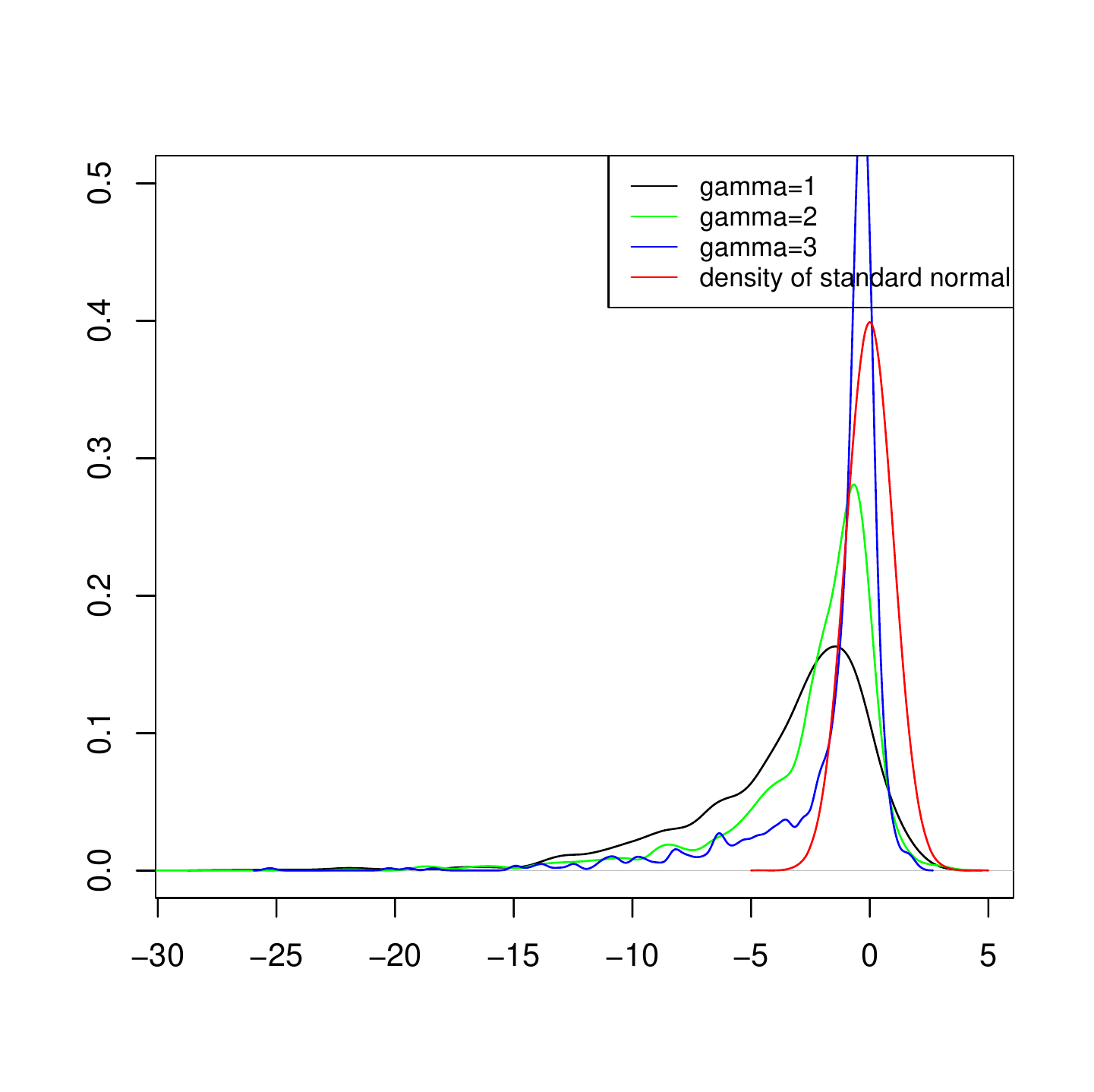}
\end{minipage}
}
\caption{\small{Asymptotic density curves of $2f(0)\mathscr{D}(\gamma)$ with $\mathscr{N}(0,1)$ innovation.}}
\end{figure}



\begin{figure}[H]
\subfloat[The near-stationary case]{
\begin{minipage}[t]{0.5\linewidth}
\centering
\includegraphics[width=6.5cm,height=6.5cm]{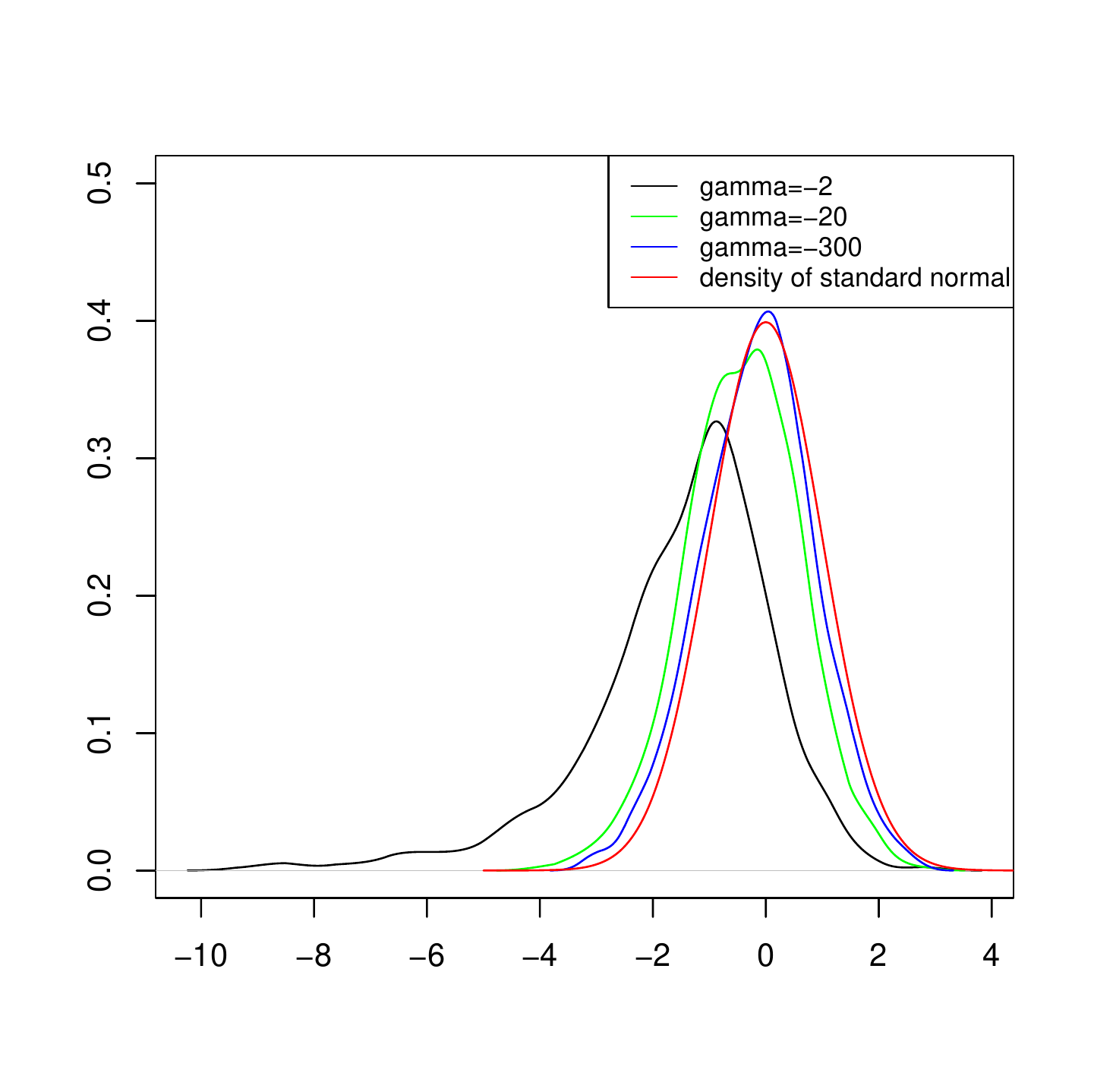}
\end{minipage}
}
\subfloat[The near-explosive case]{
\begin{minipage}[t]{0.5\linewidth}        
\includegraphics[width=6.5cm,height=6.5cm]{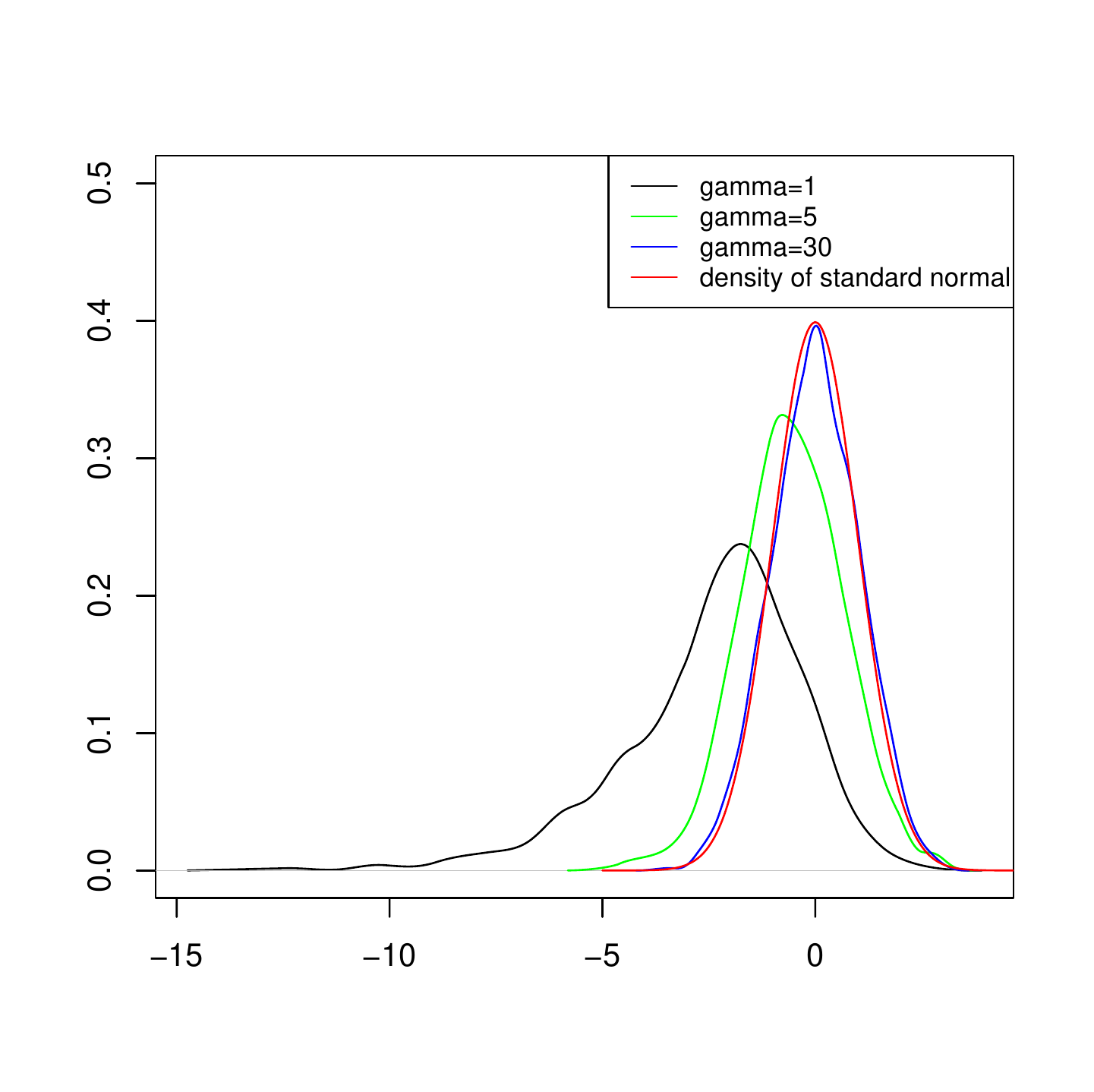}
\end{minipage}
}
\caption{\small{Asymptotic density curves of $2f(0)\mathscr{L}(\gamma)$ with $\mathscr{N}(0,1)$ innovation.}}
\end{figure}


\vskip5pt

\noindent\textbf{Accuracy of Estimators.} Andrews and Guggenberger (2008) achieved the asymptotic theory of OLS estimators for model (\ref{model}) with $n(\rho_n-1)\to 0$ and $0<\rho_n<1$. Here we will compare the accuracy between the estimators $\hat{\rho}_{\rm LAD}$ and $\hat{\rho}_{\rm OLS}$, where $\hat{\rho}_{\rm OLS}$ denotes the OLS estimator of $\rho_n$. For the near-stationary case, let $y_0=\sum_{j=0}^{\infty}\rho_n^{j}\epsilon_{-j}$, $\beta=1.1,~1.3$ and $\gamma=-5,~-50$. Here we add $10$ for each randomly selected $5\%$ sample points
to construct outliers when generating data. In these experiments, we simulate $1000$ replications each with sample size $n=200$ or $500$ to compare the empirical means (EM), absolute errors (AE) and mean squared errors (MSE) of $\hat\rho_{\rm LAD}$ and $\hat{\rho}_{\rm OLS}$ in two cases, with and without outliers. In the case without outliers, Table 1 shows that EM is very close to the true value $\rho_n$ and AE and MSE are very small, which verifies the accuracy of LAD estimator for very nearly unit root model. The errors (AE and MSE) decrease as the sample size $n$ increases. If there are outliers in the data, Table 2 shows that the errors (AE and MSE) of OLS estimator are much larger than those of LAD estimator. Furthermore, the OLS estimator occasionally produces estimates greater than one. Hence it is concluded that  $\hat{\rho}_{\rm LAD}$ is more robust than $\hat{\rho}_{\rm OLS}$.
\vskip5pt

For the near-explosive case, we take $\kappa_n=[n^{1.3}]$ in $y_0=\sum_{j=0}^{\kappa_n}\rho_n^{j}\epsilon_{-j}$, and let $\beta=1.5,~1.8$, $\gamma=5,~50$.
The selection of the parameter values here guarantees that  $n(\rho_n-1)\to 0$, $\kappa_n(1-\rho_n)\to 0$ and $\kappa_n/n\to\infty$, as $n\to \infty$.
Since the $y$ values from different cases have quite different scales, here we add $10$ times of the absolute value of the maximum for each  randomly selected $5\%$ sample points to construct outliers when generating data.
In these experiments, we simulate $1000$ replications each with sample size $n=200$ or $500$ and study the empirical means (EM), absolute errors (AE) and mean squared errors (MSE) of $\hat\rho_{\rm LAD}$ and $\hat\rho_{\rm OLS}$ in two cases, with and without outliers.
Table 3 shows that, if there are no outliers in the data, EM is very close to the true value $\rho_n$ and  AE and MSE are very small, which verifies the accuracy of the LAD estimator and OLS estimator, while, if there are outliers in the data, it is clear from Table 4 that $\hat{\rho}_{\rm LAD}$ is more robust than $\hat{\rho}_{\rm OLS}$.

\begin{table}[H]
\setlength{\abovecaptionskip}{0pt}
\setlength{\belowcaptionskip}{-0pt}
\caption{AE/EM and MSE of $\hat\rho_{\rm LAD}$($\hat\rho_{\rm OLS}$) for near-stationary AR(1) processes with $\mathscr{N}(0,1),\,U(-1,1)$ noises: no outliers case.}
\centering
\scalebox{0.85}{
\begin{tabular}{ccccccccc}
\specialrule{0.1em}{1pt}{1pt}
&&&&\multicolumn{5}{c}{$\epsilon_i$}\\
\cline{5-9}
&&&&\multicolumn{2}{c}{$N(0,1)$}&&\multicolumn{2}{c}{$U(-1,1)$}\\
\cline{5-6}\cline{8-9}
$\beta$&$\gamma$ & $n$&$\rho_n$&AE/EM&MSE&&AE/EM&MSE\\
\specialrule{0.05em}{1pt}{1pt}
1.1&-50&200&0.8528&0.0097/0.8431&0.0127&&0.0112/0.8416&0.0207\\
&&&&(0.0083/0.8445)&(0.0083)&&(0.0084/0.8445)&(0.0077)\\
&&500&0.9463&0.0037/0.9426&0.0017&&0.0049/0.9414&0.0030\\
&&&&(0.0034/0.9428)&(0.0012)&&(0.0036/0.9427)&(0.0012)\\
&-5&200&0.9853&0.0046/0.9807&0.0010&&0.0037/0.9816&0.0011\\
&&&&(0.0083/0.9770)&(0.0019)&&(0.0073/0.9780)&(0.0016)\\
&&500&0.9946&0.0026/0.9920&0.0002&&0.0025/0.9921&0.0002\\
&&&&(0.0033/0.9913)&(0.0003)&&(0.0029/0.9917)&(0.0003)\\
1.3&-50&200&0.9490&0.0072/0.9418&0.0038&&0.0101/0.9389&0.0059\\
&&&&(0.0083/0.9407)&(0.0035)&&(0.0096/0.9394)&(0.0035)\\
&&500&0.9845&0.0028/0.9817&0.0005&&0.0038/0.9807&0.0008\\
&&&&(0.0037/0.9808)&(0.0005)&&(0.0034/0.9811)&(0.0005)\\
&-5&200&0.9949&0.0029/0.9920&0.0003&&0.0025/0.9924&0.0004\\
&&&&(0.0065/0.9884)&(0.0011)&&(0.0047/0.9890)&(0.0006)\\
&&500&0.9985&0.0018/0.9966&0.0001&&0.0018/0.9967&0.0001\\
&&&&(0.0025/0.9960)&(0.0001)&&(0.0019/0.9966)&(0.0001)\\

\specialrule{0.1em}{1pt}{1pt}
\end{tabular}
}
\end{table}


\begin{table}[H]
\centering
\setlength{\abovecaptionskip}{0pt}
\setlength{\belowcaptionskip}{-0pt}
\caption{AE/EM and MSE of $\hat\rho_{\rm LAD}$($\hat\rho_{\rm OLS}$) for near-stationary AR(1) processes with $\mathscr{N}(0,1),\,U(-1,1)$ noises: outliers case.}
\scalebox{0.85}{
\begin{tabular}{ccccccccc}
\specialrule{0.1em}{1pt}{1pt}
&&&&\multicolumn{5}{c}{$\epsilon_i$}\\
\cline{5-9}
&&&&\multicolumn{2}{c}{$N(0,1)$}&&\multicolumn{2}{c}{$U(-1,1)$}\\
\cline{5-6}\cline{8-9}
$\beta$&$\gamma$ & $n$&$\rho_n$&AE/EM&MSE&&AE/EM&MSE\\
\specialrule{0.05em}{1pt}{1pt}
1.1&-50&200&0.8528&0.0046/0.8574&0.0015&&0.0046/0.8574&0.00009\\
&&&&(0.0451/0.8979)&(0.0121)&&(0.0500/0.9028)&(0.0138)\\
&&500&0.9463&0.0039/0.9501&0.0002&&0.0033/0.9495&0.0001\\
&&&&(0.0320/0.9783)&(0.0052)&&(0.0338/0.9801)&(0.0058)\\
&-5&200&0.9853&0.0015/0.9868&0.0001&&0.0013/0.9866&0.0000\\
&&&&(0.0145/0.9998)&(0.0011)&&(0.0150/1.0003)&(0.0012)\\
&&500&0.9946&0.0007/0.9954&0.0000&&0.0006/0.9953&0.0000\\
&&&&(0.0065/1.0011)&(0.0002)&&(0.0066/1.0012)&(0.0002)\\
1.3&-50&200&0.9490&0.0024/0.9514&0.0004&&0.0029/0.9519&0.0002\\
&&&&(0.0297/0.9787)&(0.0046)&&(0.0320/0.9810)&(0.0052)\\
&&500&0.9845&0.0016/0.9861&0.0000&&0.0013/0.9858&0.0000\\
&&&&(0.0137/0.9982)&(0.0009)&&(0.0140/0.9985)&(0.0010)\\
&-5&200&0.9949&0.0009/0.9958&0.0000&&0.0011/0.9960&0.0001\\
&&&&(0.0091/1.0040)&(0.0005)&&(0.0095/1.0044)&(0.0006)\\
&&500&0.9985&0.0006/0.9990&0.0001&&0.0004/0.9989&0.0001\\
&&&&(0.0040/1.0024)&(0.0001)&&(0.0040/1.0025)&(0.0001)\\

\specialrule{0.1em}{1pt}{1pt}
\end{tabular}
}
\end{table}

\begin{table}[H]
\centering
\setlength{\abovecaptionskip}{0pt}
\setlength{\belowcaptionskip}{-0pt}
\caption{AE/EM and MSE of $\hat\rho_{\rm LAD}$($\hat\rho_{\rm OLS}$) for near-explosive AR(1) processes with $\mathscr{N}(0,1),\,U(-1,1)$ noises: no outliers case.}
\scalebox{0.85}{
\begin{tabular}{ccccccccc}
\specialrule{0.1em}{1pt}{1pt}
&&&&\multicolumn{5}{c}{$\epsilon_i$}\\
\cline{5-9}
&&&&\multicolumn{2}{c}{$N(0,1)$}&&\multicolumn{2}{c}{$U(-1,1)$}\\
\cline{5-6}\cline{8-9}
$\beta$&$\gamma$ & $n$&$\rho_n$&AE/EM&MSE&&AE/EM&MSE\\
\specialrule{0.05em}{1pt}{1pt}
1.5&50&200&1.0177&2.2e-5/1.0177&2.5e-9&&2.2e-5/1.0177&2.5e-9\\
&&&&(2.2e-5/1.0177)&(2.5e-9)&&(2.2e-5/1.0177)&(2.5e-9)\\
&&500&1.0045&2.8e-5/1.0045&3.9e-9&&2.8e-5/1.0045&3.9e-9\\
&&&&(2.8e-5/1.0045)&(3.9e-9)&&(2.8e-5/1.0045)&(3.9e-9)\\
&5&200&1.0018&-6e-4/1.0012&4.5e-5&&-8e-4/1.0009&6.7e-5\\
&&&&(-e-6/1.0018)&(e-5)&&(1.8e-4/1.0019)&(e-5)\\
&&500&1.0004&-3.8e-4/1.00007&1.3e-5&&-5.7e-4/0.9999&2.8e-5\\
&&&&(1.2e-4/1.0006)&(1.5e-6)&&(0.0002/1.0007)&(e-6)\\
1.8&50&200&1.0036&-1.6e-4/1.0034&1.4e-5&&-2e-4/1.0034&1.5e-5\\
&&&&(-5e-5/1.0036)&(3e-6)&&(-e-4/1.0035)&(3.7e-6)\\
&&500&1.0007&-2.6e-4/1.0004&9e-6&&3.7e-4/1.0003&1.8e-5\\
&&&&(4e-5/1.0007)&(e-6)&&(9e-5/1.0008)&(9e-7)\\
&5&200&1.0004&-0.0017/0.9986&e-4&&-0.0020/0.9984&1.7e-4\\
&&&&(4e-4/1.0007)&(1.9e-5)&&(6e-4/1.0010)&(e-5)\\
&&500&1.00007&-9.6e-4/0.9991&3.9e-5&&-0.0011/0.9989&5.7e-5\\
&&&&(3e-4/1.00036)&(1.6e-6)&&(2.7e-4/1.00035)&(9e-7)\\

\specialrule{0.1em}{1pt}{1pt}
\end{tabular}
}
\end{table}

\begin{table}[H]
\setlength{\abovecaptionskip}{0pt}
\setlength{\belowcaptionskip}{-0pt}
\caption{AE/EM and MSE of $\hat\rho_{\rm LAD}$($\hat\rho_{\rm OLS}$) for near-explosive AR(1) processes with $\mathscr{N}(0,1),\,U(-1,1)$ noises: outliers case.}
\centering
\scalebox{0.85}{
\begin{tabular}{ccccccccc}
\specialrule{0.1em}{1pt}{1pt}
&&&&\multicolumn{5}{c}{$\epsilon_i$}\\
\cline{5-9}
&&&&\multicolumn{2}{c}{$N(0,1)$}&&\multicolumn{2}{c}{$U(-1,1)$}\\
\cline{5-6}\cline{8-9}
$\beta$&$\gamma$ & $n$&$\rho_n$&AE/EM&MSE&&AE/EM&MSE\\
\specialrule{0.05em}{1pt}{1pt}
1.5&50&200&1.0177&-0.0144/1.0033&0.0020&&-0.0140/1.0037&0.0019\\
&&&&(-0.4835/0.5341)&(2.1195)&&(-0.4703/0.5474)&(2.0560)\\
&&500&1.0045&-7.1e-3/0.9974&5.1e-4&&-0.0069/0.9976&5e-4\\
&&&&(-0.5434/0.4610)&(1.9600)&&(-0.5345/0.4700)&(1.9105)\\
&5&200&1.0018&-0.0080/0.9939&3.8e-4&&-0.0090/0.9928&4.6e-4\\
&&&&(-0.7847/0.2171)&(3.1892)&&(-0.7748/0.2270)&(3.1400)\\
&&500&1.0004&-0.0064/0.9941&2.5e-4&&-0.0074/0.9930&3.2e-4\\
&&&&(-0.7823/0.2182)&(3.1841)&&(-0.7731/0.2273)&(3.1062)\\
1.8&50&200&1.0036&-0.0033/1.0003&1.2e-4&&-0.0038/0.9998&1.5e-4\\
&&&&(-0.7888/0.2148)&(3.1625)&&(-0.7929/0.2108)&(3.1825)\\
&&500&1.0007&-0.0045/0.9962&1.5e-4&&-0.0055/0.9952&2.1e-4\\
&&&&(-0.7930/0.2077)&(3.2306)&&(-0.7952/0.2055)&(3.2533)\\
&5&200&1.0004&-9.6e-3/0.9907&4.8e-4&&-9.8e-3/0.9905&5e-4\\
&&&&(-0.7404/0.2600)&(2.9465)&&(-0.7364/0.2640)&(2.9257)\\
&&500&1.00007&-8.4e-3/0.9917&3.8e-4&&-9.0e-3/0.9911&4.3e-4\\
&&&&(-0.7375/0.2625)&(2.9256)&&(-0.7394/0.2607)&(2.9476)\\

\specialrule{0.1em}{1pt}{1pt}
\end{tabular}
}
\end{table}



\section{Concluding remarks}\label{s5}

In this article we develop the limit theory of the LAD estimator for the AR(1) process with a root close to unity. The parameter $\rho_n$ satisfies $n(\rho_n-1)\to 0$ or $n(\rho_n-1)\to \gamma$ for some fixed number $\gamma\neq 0$, as $n\to \infty$. It is
shown that, in the first case, for both near-stationary and near-explosive processes, the
LAD estimator and the $t$-type statistic have Cauchy and normal asymptotic distribution, respectively.
The simulation study in this case confirms the theoretical results and it illustrates that the theory should be useful in statistical inference and the LAD estimator is robust if there are outliers in the data. In the case that $n(\rho_n-1)\to \gamma$ for some fixed number $\gamma\ne 0$, as $n\to \infty$, we also develop asymptotic theory for the LAD estimator and the $t$-type statistic under the assumption that $y_0=0$  which gives us the connection with the existing results in literature, e.g. Chan and Wei (1987), Phillips (1987), and Herce (1996).
\vskip5pt

In summary, we establish the asymptotic theory of the LAD estimators for the nearly unit root processes which correspond to the results of the OLS estimators developed by Chan and Wei (1987), Phillips (1987, 1988, 2021), Phillips and Magdalinos (2007), Andrews and Guggenberger (2008), Chan (2009) and so on.
The results in this article are also the completion of the LAD estimators studied in Zhou and Lin (2014) and Wang {\it et al.} (2020). These authors worked on the case that $\rho_n\to1$ and $n|\rho_n-1|\to\infty$ as $n\to\infty$.  The normalizer of the  asymptotic distributions for the LAD estimators there has a rate in the range $(n^{1/2}, n)$ for the near-stationary case. We enlarge the rate spectrum of the normalizer in the asymptotic distributions of the LAD estimators.



\section{Proof of main results}\label{s6}


We first give some comments on the proofs. It is known that the methods of the LAD estimation are classic which were developed by Pollard (1991), Davis {\it et al.} (1992), Knight (1989, 1998) and Ling (2005). The limiting behaviors of the quadratic functionals, $\sum_{i=1}^ny_{i-1}^2$ and $\sum_{i=1}^ny_{i-1}{\rm sign}(\epsilon_i)$, play the crucial role in the analysis. Concretely, for the very nearly case, $n(\rho_n-1)\to0$, we mainly follow the strategies of Zhou and Lin (2014) and Wang {\it et al.} (2020), while, under our framework, the initial values $y_0$ dominate the asymptotic behavior of $\sum_{i=1}^{n}y_{i-1}^2$; see Lemma \ref{l-s0} and Lemma \ref{l-e0}. More work is required for the local unit root case, $n(\rho_n-1)\to\gamma$ for some $\gamma\neq0$. To treat the asymptotic joint distributions of $\left(\frac{1}{n}\sum_{i=1}^ny_{i-1}{\rm sign}(\epsilon_i),\frac{1}{n^2}\sum_{i=1}^ny_{i-1}^2\right)$, we need to develop a functional central limit theorem for the process $\{(K_n(t),L_n(t)),0\leq t\leq1\}$ defined as in (\ref{kn}), i.e. Lemma \ref{KLM1}, by the vector-value martingale invariance principle. We remark that, although Phillips and Durlauf (1986) (Lemma 3.1) established the asymptotic theory for sample moments of vector-value integrated processes, their results can not be used here because the regression parameter $\rho_n$ depends on the sample size $n$ which causes the covariance matrix defined as in (\ref{covm}) to be time-dependent. Finally, by using the tools from stochastic calculus, we achieve the estimates, Lemma \ref{result3-lem1} and Lemma \ref{result3-lem2}, which is the key of the proof in Theorem \ref{main result3}.
\vskip5pt

Throughout the proof, we use the following identity of Knight (1998),
\begin{align}\label{Knight}
{\rm for}~x\neq0, \quad |x-y|-|x|=-y{\rm sign}(x)+2\int_{0}^{y}\left(\mathbb{I}_{\{x\leq s\}}-\mathbb{I}_{\{x\leq 0\}}\right){\rm d}s.
\end{align}
From the model (\ref{model}), we also observe that
\begin{align}\label{near-s-model}
y_{i}=\rho_n^{i}y_0+\sum_{j=1}^{i}\rho_n^{i-j}\epsilon_j, \quad 1\leq i\leq n.
\end{align}


\subsection{The case $n(\rho_n-1)\to0$ as $n\to\infty$}
In this subsection, we will prove Theorem \ref{near-s0} and Theorem \ref{near-e0}.
Let us start by presenting two lemmas proved in Appendix.

\begin{lemma}\label{l-s0}
Under the assumptions of Theorem \ref{near-s0}, we have
\begin{itemize}
\item[(1).] $\sqrt{1-\rho_n^2}y_0\rightarrow_d\sigma X$;
\item[(2).] $\frac{1}{\sqrt{n}}\sum_{i=1}^n\rho_n^{i-1}{\rm sign}(\epsilon_i)\rightarrow_d Y$;
\item[(3).] $\sqrt{(1-\rho_n^2)/n}\max_{1\leq i\leq n}|y_{i-1}|\rightarrow_p0$;
\end{itemize}
here $X$ and $Y$ are independent standard normal random variables.
\end{lemma}

\begin{lemma}\label{l-e0}
Under the assumptions of Theorem \ref{near-e0}, we have
\begin{itemize}
\item[(1).] $\frac{1}{\sqrt{\kappa_n}}y_0\rightarrow_d\sigma U$;
\item[(2).] $\frac{1}{\sqrt{n}}\sum_{i=1}^n\rho_n^{i-1}{\rm sign}(\epsilon_i)\rightarrow_d V$;
\item[(3).] $\frac{1}{\sqrt{n\kappa_n}}\max_{1\leq i\leq n}|y_{i-1}|\rightarrow_p 0$;
\end{itemize}
here $U$ and $V$ are independent standard normal random variables.
\end{lemma}

\begin{proof}[\textbf{Proof of Theorem \ref{near-s0}}]
Denote $\hat{u}_n=\sqrt{n/(1-\rho_n^2)}(\hat{\rho}_{\rm LAD}-\rho_n)$. Then $\hat{u}_n$ is the minimizer of the following convex objective function,
\begin{align*}
Z_n(u)=\sum_{i=1}^n\left(|\epsilon_i-\sqrt{(1-\rho_n^2)/n}y_{i-1}u|-|\epsilon_i|\right).
\end{align*}
Based on the ideas of Davis {\it et al.} (1992) and Ling (2005), if we can prove that, for each $u$, $Z_n(u)$ converges weakly to a random variable which has a unique minimizer $u_{\min}$, then $\hat{u}_n$ must converge weakly to $u_{\min}$.
\vskip5pt

By (\ref{Knight}), we rewrite $Z_n(u)$ as follows
\begin{align}\label{decom}
Z_n(u)=-uA_n+2\sum_{i=1}^n\xi_i(u),
\end{align}
here
\[
A_n=\sqrt{(1-\rho_n^2)/n}\sum_{i=1}^ny_{i-1}{\rm sign}(\epsilon_i)
\]
and
\[
\xi_i(u)=\int_0^{\sqrt{(1-\rho_n^2)/n}y_{i-1}u}\left(\mathbb{I}_{\{\epsilon_i\leq s\}}-\mathbb{I}_{\{\epsilon_i\leq 0\}}\right){\rm d}s.
\]
Next we analyze the asymptotic properties of $A_n$ and $\xi_i(u)$, respectively.
By (\ref{near-s-model}), we can decompose $A_n$ as follows
\begin{align}\label{An0}
A_n&=\sqrt{(1-\rho_n^2)/n}\sum_{i=1}^{n}\Big(\rho_n^{i-1}y_0+\sum_{j=1}^{i-1}\rho_n^{i-1-j}\epsilon_j\Big)\,{\rm sign}(\epsilon_i)\nonumber\\
&=\sqrt{(1-\rho_n^2)}y_0\cdot\frac{1}{\sqrt{n}}\sum_{i=1}^n\rho_n^{i-1}{\rm sign}(\epsilon_i)+R_{n,1},
\end{align}
where the remainder
\[
R_{n,1}=\sqrt{(1-\rho_n^2)/n}\sum_{i=1}^n\sum_{j=1}^{i-1}\rho_n^{i-1-j}\epsilon_j\,{\rm sign}(\epsilon_i).
\]
The parts (1) and (2) of Lemma \ref{l-s0}, together with the continuous mapping theorem, yield that the first term in the equation (\ref{An0}) converges weakly to the random variable $\sigma XY$, where $X$ and $Y$ are independent standard normal random variables. Now, we show that the remainder converges to zero in probability. Obviously, by simple calculations, $\mathbb{E}(R_{n,1})=0$ and
\begin{align*}
\mathbb{E}(R_{n,1}^2)&=\frac{1-\rho_n^2}{n}\,\mathbb{E}\Big[\Big(\sum_{i=1}^{n}\sum_{j=1}^{i-1}\rho_n^{i-1-j}\epsilon_j\,{\rm sign}(\epsilon_i)\Big)^2\Big]\\
&=\frac{1-\rho_n^2}{n}\sum_{i=1}^{n}\mathbb{E}\Big[\Big(\sum_{j=1}^{i-1}\rho_n^{i-1-j}\epsilon_j\Big)^2\Big]
=\frac{n-1}{n}\Big(1-\frac{\rho_n^2-\rho_n^{2n}}{(n-1)(1-\rho_n^2)}\Big)\sigma^2\to0
\end{align*}
as $n\to\infty$, where the last step is due to Proposition A.1. Consequently, $A_n$ converges weakly to the random variable $\sigma XY$.
\vskip5pt

We now analyze the second term $\sum_{i=1}^n\xi_i(u)$ in the equation (\ref{decom}) by the martingale method. Denote the filtration by $\mathcal{F}_{n,i}=\sigma(y_0,\epsilon_1,\ldots,\epsilon_i)$ for $1\leq i\leq n$ and $n\geq1$, then we can write
\begin{align*}
\sum_{i=1}^n\xi_i(u)=\sum_{i=1}^n\mathbb{E}(\xi_i(u)|\mathcal{F}_{n,i-1})+R_{n,2},
\end{align*}
where the remainder
\[
R_{n,2}=\sum_{i=1}^n(\xi_i(u)-\mathbb{E}(\xi_i(u)|\mathcal{F}_{n,i-1})).
\]
is a martingale with respect to the filtration, $\{\mathcal{F}_{n,i},1\leq i\leq n\}$.
\vskip5pt

By the Taylor's formula,
\begin{align}\label{xi1}
\sum_{i=1}^n\mathbb{E}(\xi_i(u)|\mathcal{F}_{n,i-1})&=\sum_{i=1}^n\int_0^{\sqrt{(1-\rho_n^2)/n}y_{i-1}u}\mathbb{E}\left(\mathbb{I}_{\{\epsilon_i\leq s\}}-\mathbb{I}_{\{\epsilon_i\leq 0\}}\right){\rm d}s\nonumber\\
&=\sum_{i=1}^n\int_0^{\sqrt{(1-\rho_n^2)/n}y_{i-1}u}\left(F(s)-F(0)\right){\rm d}s\nonumber\\
&=\sum_{i=1}^n\int_0^{\sqrt{(1-\rho_n^2)/n}y_{i-1}u}\Big(f(0)s+\frac{1}{2}f'(s^*)s^2\Big){\rm d}s\nonumber\\
&=\frac{u^2f(0)}{2}\cdot\frac{1-\rho_n^2}{n}\sum_{i=1}^ny_{i-1}^2+R_{n,3},
\end{align}
where $F(\cdot)$ is the distribution function of $\epsilon_1$, $s^*\in(0,s)$ and the remainder
\[
R_{n,3}=\frac{1}{2}\sum_{i=1}^n\int_0^{\sqrt{(1-\rho_n^2)/n}y_{i-1}u} f'(s^*)s^2 {\rm d}s.
\]

We first consider the sample variance, $\sum_{i=1}^ny_{i-1}^2$. From equation (\ref{near-s-model}),
\begin{align}\label{y_t2}
\sum_{i=1}^ny_{i-1}^2=\sum_{i=1}^n\rho_n^{2(i-1)}y_0^2+\sum_{i=1}^n\rho_n^{i-1}\sum_{j=1}^{i-1}\rho_{n}^{i-j-1}\epsilon_jy_0
+\sum_{i=1}^n\Big(\sum_{j=1}^{i-1}\rho_{n}^{i-j-1}\epsilon_j\Big)^2.
\end{align}
Simple calculations yield that
\begin{align*}
\mathbb{E}\Big[\sum_{i=1}^n\Big(\sum_{j=1}^{i-1}\rho_{n}^{i-j-1}\epsilon_j\Big)^2\Big]=\sum_{i=1}^n\sum_{j=1}^{i-1}\rho_n^{2(i-1-j)}\sigma^2
=\frac{(n-1)\sigma^2}{1-\rho_n^2}\Big(1-\frac{\rho_n^2-\rho_n^{2n}}{(n-1)(1-\rho_n^2)}\Big).
\end{align*}
Note that, by the Cauchy-Schwarz's inequality,
\[
\Big|\sum_{i=1}^n\rho_n^{i-1}\sum_{j=1}^{i-1}\rho_{n}^{i-j-1}\epsilon_j\Big|^2\leq\sum_{i=1}^n\rho_n^{2(i-1)}
\sum_{i=1}^n\Big(\sum_{j=1}^{i-1}\rho_{n}^{i-j-1}\epsilon_j\Big)^2,
\]
so, applying Proposition A.1 for $k_n=n$ and part (1) of Lemma \ref{l-s0}, for the second term in (\ref{y_t2}), we can obtain
\[
\frac{1-\rho_n^2}{n}\sum_{i=1}^n\rho_n^{i-1}\sum_{j=1}^{i-1}\rho_{n}^{i-j-1}\epsilon_j y_0\rightarrow_p0.
\]
It also follows, for the third term in (\ref{y_t2}),
\[
\frac{1-\rho_n^2}{n}\sum_{i=1}^n\Big(\sum_{j=1}^{i-1}\rho_{n}^{i-j-1}\epsilon_j\Big)^2\rightarrow_p0.
\]
Moreover, Proposition A.1 and part (1) of Lemma \ref{l-s0}, combined with the continuous mapping theorem, imply that
\[
\frac{1-\rho_n^2}{n}\sum_{i=1}^n\rho_n^{2(i-1)}y_0^2\rightarrow_d\sigma^2X^2.
\]
Therefore, we have shown that
\begin{align}\label{sample-var0}
\frac{1-\rho_n^2}{n}\,\sum_{i=1}^ny_{i-1}^2\rightarrow_d\sigma^2X^2.
\end{align}
\vskip5pt

It remains to deal with the remainders, $R_{n,2}$ and $R_{n,3}$. Notice that
\begin{align*}
|R_{n,3}|&\leq\frac{|u|^3\sup_{x\in\mathbb{R}}|f'(x)|}{3}\left(\frac{1-\rho_n^2}{n}\right)^{3/2}\sum_{i=1}^n|y_{i-1}|^3\\
&\leq\frac{|u|^3\sup_{x\in\mathbb{R}}|f'(x)|}{3}\left(\frac{1-\rho_n^2}{n}\right)^{3/2}\max_{1\leq i\leq n}|y_{i-1}|\sum_{i=1}^ny_{i-1}^2,
\end{align*}
by equation (\ref{sample-var0}) and part (3) of Lemma \ref{l-s0}, $R_{n,3}$ converges to zero in probability. Since $R_{n,2}$ is a martingale, with the methods of Wang {\it et al.} (2020) or Zhou and Lin (2014) and the help of Lemma \ref{l-s0}, we can prove that $R_{n,2}$ converges to zero in probability.
\vskip5pt

Summarizing above discussions, we get that
\begin{align*}
Z_n(u)&=-uA_n+2\sum_{i=1}^n\xi_i(u)\nonumber\\
&=-u\sqrt{(1-\rho_n^2)}y_0\cdot\frac{1}{\sqrt{n}}\sum_{i=1}^n\rho_n^{i-1}{\rm sign}(\epsilon_i)
+u^2f(0)\frac{1-\rho_n^2}{n}\sum_{i=1}^n\rho_n^{2(i-1)}y_0^2+o_p(1)\nonumber\\
&\rightarrow_d-u \sigma XY+u^2\sigma^2f(0)X^2=:Z(u).
\end{align*}
Because $Z(u)$ has a unique minimum at
\[
u_{\min}=Y/(2\sigma f(0)X),
\]
by Lemma 2.2 of Davis {\it et al.} (1992), we derive that
\begin{align*}
\hat{u}_n=\sqrt{n/(1-\rho_n^2)}(\hat{\rho}_{\rm LAD}-\rho_n)\rightarrow_d\frac{1}{2\sigma f(0)}\mathcal{C}.
\end{align*}
Combining this with equation (\ref{sample-var0}) and using the continuous mapping theorem, we complete the proof of Theorem \ref{near-s0}.
\end{proof}

\begin{proof}[\textbf{Proof of Theorem \ref{near-e0}}]
The proof is similar to that of Theorem \ref{near-s0}. In this case,  $\hat{u}_n=\sqrt{n\kappa_n}(\hat{\rho}_{\rm LAD}-\rho_n)$,
$Z_n(u)=\sum_{i=1}^n\left(|\epsilon_i-(n\kappa_n)^{-1/2}y_{i-1}u|-|\epsilon_i|\right)$, and in (\ref{decom}),
\[
A_n=\frac{1}{\sqrt{n\kappa_n}}\sum_{i=1}^ny_{i-1}{\rm sign}(\epsilon_i)
\]
and
\[
\xi_i(u)=\int_0^{y_{i-1}u/\sqrt{n\kappa_n}}\left(\mathbb{I}_{\{\epsilon_i\leq s\}}-\mathbb{I}_{\{\epsilon_i\leq 0\}}\right){\rm d}s.
\]
Lemma \ref{l-e0} is applied to prove the corresponding terms as in the proof of Theorem \ref{near-s0}.  Here we only partially demonstrate the role of the assumptions, $\kappa_n/n\to\infty$ and $\kappa_n(\rho_n-1)\to0$ as $n\to\infty$. Let us consider the normalized sample variance. For the third term of (\ref{y_t2}), we want to show that
\[
\frac{1}{n\kappa_n}\sum_{i=1}^n\Big(\sum_{j=1}^{i-1}\rho_{n}^{i-j-1}\epsilon_j\Big)^2\rightarrow_p0.
\]
By a simple calculation, we get
\begin{align*}
\mathbb{E}\Big[\sum_{i=1}^n\Big(\sum_{j=1}^{i-1}\rho_{n}^{i-j-1}\epsilon_j\Big)^2\Big]
=\sum_{i=1}^n\sum_{j=1}^{i-1}\rho_n^{2(i-1-j)}\sigma^2
=\frac{(n-1)\sigma^2}{1-\rho_n^2}\Big(1-\frac{\rho_n^2-\rho_n^{2n}}{(n-1)(1-\rho_n^2)}\Big).
\end{align*}
Then, Proposition A.1 yields that
\begin{align*}
\frac{1}{n\kappa_n}\cdot\frac{(n-1)\sigma^2}{1-\rho_n^2}\Big(1-\frac{\rho_n^2-\rho_n^{2n}}{(n-1)(1-\rho_n^2)}\Big)
=O\Big(\frac{n}{\kappa_n}\Big),
\end{align*}
which tends to zero if $\kappa_n/n\to\infty$ as $n\to\infty$. The proof is complete.
\end{proof}


\subsection{The case $n(\rho_n-1)\to\gamma$ as $n\to\infty$}

In this subsection, we will prove Theorem \ref{main result2} and Theorem \ref{main result3}.
\vskip5pt

Firstly, we prove Theorem \ref{main result2} following the line as in the proofs of the case $n(\rho_n-1)\to 0$.
In this case, $\hat{u}_n=n(\hat{\rho}_{\rm LAD}-\rho_n)$ and it is the minimizer of the following convex function,
\begin{align}\label{zn2}
Z_n(u)=\sum_{i=1}^n\left(|\epsilon_i-n^{-1}y_{i-1}u|-|\epsilon_i|\right).
\end{align}
Again by the Knight's identity (\ref{Knight}), we can rewrite the quantity $Z_n(u)$ in (\ref{zn2}) as identity (\ref{decom}) with
\begin{align}\label{A_n}
A_n=\frac{1}{n}\sum_{i=1}^ny_{i-1}{\rm sign}(\epsilon_i)
\end{align}
and
\begin{align}
\xi_i(u)=\int_0^{n^{-1}y_{i-1}u}\left(\mathbb{I}_{\{\epsilon_i\leq s\}}-\mathbb{I}_{\{\epsilon_i\leq 0\}}\right){\rm d}s.
\end{align}

Again we use the filtration $\mathcal{F}_{n,i}=\sigma(y_0,\epsilon_1,\ldots,\epsilon_i)$ for $1\leq i\leq n$ and $n\geq1$, and write
\begin{align*}
\sum_{i=1}^n\xi_i(u)=\sum_{i=1}^n\mathbb{E}(\xi_i(u)|\mathcal{F}_{n,i-1})+R_{n,1},
\end{align*}
where
\begin{align}\label{rn2}
R_{n,1}=\sum_{i=1}^n(\xi_i(u)-\mathbb{E}(\xi_i(u)|\mathcal{F}_{n,i-1}))
\end{align}
is a martingale with respect to the filtration, $\{\mathcal{F}_{n,i},1\leq i\leq n\}$.
By the same argument as in (\ref{xi1}),
\begin{align*}
\sum_{i=1}^n\mathbb{E}(\xi_i(u)|\mathcal{F}_{n,i-1})
&=\frac{u^2f(0)}{2}\cdot\frac{1}{n^2}\sum_{i=1}^ny_{i-1}^2+R_{n,2},
\end{align*}
where
\begin{align}\label{rn3}
R_{n,2}=\frac{1}{2}\sum_{i=1}^n\int_0^{n^{-1}y_{i-1}u} f'(s^*)s^2 {\rm d}s
\end{align}
and $s^*$ is between $0$ and $s$.
As in the proofs of Theorem \ref{near-s0} and Theorem \ref{near-e0}, to obtain the asymptotic distribution of $Z_n(u)$, the key point is to analyze the joint distribution of $\left(\frac{1}{n}\sum_{i=1}^ny_{i-1}{\rm sign}(\epsilon_i),\frac{1}{n^2}\sum_{i=1}^ny_{i-1}^2\right)$.

\begin{lemma}\label{KLM1}
Under the  assumptions of Theorem \ref{main result2},
the process $(K_n(t),L_n(t), 0\leq t\leq1)$ defined as in (\ref{kn}) converges weakly to the continuous process
 $({\bf{X}}(t), 0\leq t\leq1)$ defined as in (\ref{xt})
with independent Gaussian increments, mean vector zero and covariance matrix $\Gamma(t)$ defined as in (\ref{covm}). Moreover, we also have
\begin{align}\label{KLM1-0}
{\bf X}^{\rm T}(t)=\int_0^t\Lambda^{1/2}(s)\,{\rm d}{\bf B}^{\rm T}(h(s)),\quad 0\leq t\leq1,
\end{align}
where $\Lambda(\cdot)$ is a non-negative definite matrix-valued function given by
\begin{align*}
\Lambda(t)=(\lambda_{ij}(t)):=\frac{1}{2}\left(
             \begin{array}{cc}
               {1-\tanh(2\gamma(1-t))}  & {\mathbb{E}|\,\epsilon_1|}{{\rm sech}(2\gamma(1-t))} \\
               {\mathbb{E}|\,\epsilon_1|}{{\rm sech}(2\gamma(1-t))} & {1+\tanh(2\gamma(1-t))} \\
             \end{array}
           \right), \quad 0\leq t\leq1,
\end{align*}
${\bf B}=(B_1,B_2)$ is a $2$-dimensional standard Brownian motion and
\begin{align*}
h(t)=\frac{1}{\gamma}\big(\sinh(2\gamma)-\sinh(2\gamma(1-t))\big), \quad 0\leq t\leq1.
\end{align*}
\end{lemma}

\begin{proposition}\label{joint conv2}
Under the assumptions of Theorem \ref{main result2}, as $n\to\infty$,
\begin{align*}
\left(\frac{1}{n}\sum_{i=1}^ny_{i-1}{\rm sign}(\epsilon_i),\frac{1}{n^2}\sum_{i=1}^ny_{i-1}^2\right)\rightarrow_d
\left(\int_0^1L(t)\,{\rm d}K(t), \int_0^1e^{-2\gamma(1-t)}L^2(t)\,{\rm d}t\right).
\end{align*}
\end{proposition}

\begin{proof}
Observe that
\begin{align*}
\frac{1}{n}\sum_{i=1}^ny_{i-1}{\rm sign}(\epsilon_i)&=\sum_{i=1}^n\left(\frac{1}{\sqrt{n}}\sum_{j=1}^{i-1}\rho_n^{n-j}\epsilon_j\right)\frac{1}{\sqrt{n}}\rho_n^{i-1-n}{\rm sign}(\epsilon_i)\\
&=\int_0^1L_n(t)\,{\rm d}K_n(t)
\end{align*}
and
\begin{align*}
\frac{1}{n^2}\sum_{i=1}^ny_{i-1}^2&=\sum_{i=1}^n\frac{1}{n}\left(\frac{1}{\sqrt{n}}\sum_{j=1}^i\rho_n^{i-j}\epsilon_j\right)^2\\
&=\int_0^1e^{-2\gamma(1-t)}L_n^2(t)\,{\rm d}t+R_{n,3},
\end{align*}
where the remainder
\[
R_{n,3}=\frac{1}{n}\sum_{i=1}^n\rho_n^{2(i-n)}L_n^2\left(\frac{i}{n}\right)-\int_0^1e^{-2\gamma(1-t)}L_n^2(t)\,{\rm d}t.
\]
According to the same argument in the proof of Lemma 2.2 in Chan and Wei (1987), we can show that $R_{n,3}\rightarrow_p0$ as $n\to\infty$. Hence Theorem 2.1 in Hansen (1992) and Lemma \ref{KLM1} immediately yield the desired result.
\end{proof}

\begin{proof}[\textbf{Proof of Theorem \ref{main result2}}]
Recall the analysis in the beginning of this subsection. Notice that $R_{n,2}$ in (\ref{rn3}) satisfies
\begin{align}\label{2rn2}
2|R_{n,2}|&\leq\frac{|u|^3\sup_{x\in\mathbb{R}}|f'(x)|}{3}\frac{1}{n^3}\sum_{i=1}^n|y_{i-1}|^3\notag\\
&\leq\frac{|u|^3\sup_{x\in\mathbb{R}}|f'(x)|}{3}\frac{1}{n}\max_{1\leq i\leq n}|y_{i-1}|\cdot\frac{1}{n^2}\sum_{i=1}^ny_{i-1}^2.
\end{align}
By Proposition \ref{joint conv2} and Proposition A.2, $R_{n,2}$ converges to zero in probability.
In addition, following the same line of Wang {\it et al.} (2020) or Zhou and Lin (2014), we can also show that $R_{n,1}$ in (\ref{rn2}) converges to zero in probability. Combining  with Proposition \ref{joint conv2}, we have
\begin{align*}
Z_n(u)&=-uA_n+2\sum_{i=1}^n\xi_i(u)\nonumber\\
&\rightarrow_d -u\int_0^1L(t)\,{\rm d}K(t)+{u^2f(0)}\int_0^1e^{-2\gamma(1-t)}L^2(t)\,{\rm d}t=:Z(u).
\end{align*}
Because $Z(u)$ has a unique minimum at
\[
u_{\min}=\frac{1}{2f(0)}\cdot\frac{\int_0^1L(t)\,{\rm d}K(t)}{\int_0^1e^{-2\gamma(1-t)}L^2(t)\,{\rm d}t},
\]
by Lemma 2.2 of Davis {\it et al.} (1992),
\begin{align*}
\hat{u}_n=n(\hat{\rho}_{\rm LAD}-\rho_n)\rightarrow_d\frac{1}{2f(0)}\cdot\frac{\int_0^1L(t)\,{\rm d}K(t)}{\int_0^1e^{-2\gamma(1-t)}L^2(t)\,{\rm d}t},
\end{align*}
as $n\to\infty$. Applying the continuous mapping theorem and Proposition \ref{joint conv2} again,
we complete the proof of Theorem \ref{main result2}.
\end{proof}


\begin{proof}[\textbf{Proof of Theorem \ref{main result3}}]
Denote the square root of the matrix-valued function $\Lambda(t)$ defined as in Lemma \ref{KLM1} by
\[
\Lambda^{1/2}(t)=\left(
                   \begin{array}{cc}
                     \widetilde{\lambda}_{11}(t) & \widetilde{\lambda}_{12}(t) \\
                     \widetilde{\lambda}_{12}(t) & \widetilde{\lambda}_{22}(t) \\
                   \end{array}
                 \right), \quad 0\leq t\leq1.
\]
Then, from Lemma \ref{KLM1}, we know that
\begin{align}\label{main result3-0}
\left\{
\begin{array}{cc}
  &{\rm d}K(t)=\widetilde{\lambda}_{11} (t){\rm d}B_1(h(t))+\widetilde{\lambda}_{12} (t){\rm d}B_2(h(t))\\
  &{\rm d}L(t)=\widetilde{\lambda}_{12} (t){\rm d}B_1(h(t))+\widetilde{\lambda}_{22} (t){\rm d}B_2(h(t))
\end{array}
\right..
\end{align}
Moreover, by the time change for It\^{o} integrals, e.g. Theorem 8.5.7 in {\O}ksendal (2005), we can rewrite equations (\ref{main result3-0}) as
\begin{align}\label{main result3-1}
\left\{
\begin{array}{cc}
  &{\rm d}K(t)=\widetilde{\lambda}_{11} (t)\sqrt{h'(t)}{\rm d}\widetilde{B}_1(t)+\widetilde{\lambda}_{12} (t)\sqrt{h'(t)}{\rm d}\widetilde{B}_2(t)\\
  &{\rm d}L(t)=\widetilde{\lambda}_{12} (t)\sqrt{h'(t)}{\rm d}\widetilde{B}_1(t)+\widetilde{\lambda}_{22} (t)\sqrt{h'(t)}{\rm d}\widetilde{B}_2(t)
\end{array}
\right.,
\end{align}
here $\widetilde{B}=(\widetilde{B}_1,\widetilde{B}_2)$ is also a 2-dimensional standard Brownian motion. If denote
\begin{align*}
{\rm d}\widehat{B}(t)=\frac{\widetilde{\lambda}_{11} (t)\sqrt{h'(t)}{\rm d}\widetilde{B}_1(t)+\widetilde{\lambda}_{12} (t)\sqrt{h'(t)}{\rm d}\widetilde{B}_2(t)}{\sqrt{\widetilde{\lambda}_{11}^2(t)h'(t)+\widetilde{\lambda}_{12}^2(t)h'(t)}},
\end{align*}
then, by the L\'{e}vy characterization of Brownian motion, e.g. Theorem 8.6.1 in {\O}ksendal (2005), we know that $\widehat{B}$ is a 1-dimensional standard Brownian motion. Notice that
\[
\widetilde{\lambda}_{11}^2(t)h'(t)+\widetilde{\lambda}_{12}^2(t)h'(t)=e^{-2\gamma(1-t)}, \quad 0\leq t\leq1,
\]
hence the random variable $\mathscr{L}(\gamma)$ can be represented as
\begin{align*}
2f(0)\cdot\mathscr{L}(\gamma)=\frac{\int_0^1e^{-\gamma(1-t)}L(t)\,{\rm d}\widehat{B}(t)}{\sqrt{\int_0^1e^{-2\gamma(1-t)}L^2(t)\,{\rm d}t}}.
\end{align*}

Now, the subsequent proof will base on Theorem 1 of Rootz\'{e}n (1980) which is restated as Theorem A.2 in Appendix. We need verify the assumptions in Theorem A.2. Let
\begin{align*}
\varphi_\gamma(t)=2\gamma e^{-\gamma}e^{-\gamma(1-t)}L(t), \quad 0\leq t\leq1, \; \gamma>0,
\end{align*}
and
\begin{align*}
\psi_\gamma(t)=\sqrt{-2\gamma}e^{-\gamma(1-t)}L(t), \quad 0\leq t\leq1, \; \gamma<0,
\end{align*}
then we have the following lemmas whose proofs are postponed to Appendix.

\begin{lemma}\label{result3-lem1}
As $\gamma\to\infty$,
\begin{align}\label{result3-lem1-1}
\sup_{0\leq t\leq 1}\left|\int_0^t\varphi_\gamma(s)\,{\rm d}s\right|\rightarrow_p0
\end{align}
and
\begin{align}\label{result3-lem1-2}
\int_0^1\varphi_\gamma^2(t)\,{\rm d}t\rightarrow_p \mathscr{N}^2(0,1).
\end{align}
\end{lemma}

\begin{lemma}\label{result3-lem2}
As $\gamma\to-\infty$,
\begin{align}\label{result3-lem2-1}
\sup_{0\leq t\leq 1}\left|\int_0^t\psi_\gamma(s)\,{\rm d}s\right|\rightarrow_p0
\end{align}
and
\begin{align}\label{result3-lem2-2}
\int_0^1\psi_\gamma^2(t)\,{\rm d}t\rightarrow_p 1.
\end{align}
\end{lemma}

With the help of Lemma \ref{result3-lem1}, Lemma \ref{result3-lem2}, and Theorem A.2, we complete the proof of Theorem \ref{main result3}.
\end{proof}



\section*{Technical appendix and proofs}\label{s7}

\noindent\textbf{Proposition A.1.}\label{app pro} {\it Let $\{k_n,n\geq1\}$ be a sequence of positive numbers such that $k_n\to\infty$ and $k_n|1-\rho_n^2|\to0$ as $n\to\infty$, then we have $\rho_n^{2k_n}\to1$ and
\begin{align}\label{app0}
\frac{1-\rho_n^{2k_n}}{k_n(1-\rho_n^2)}=1+O(k_n(1-\rho_n^2))
\end{align}
as $n\to\infty$.}

\begin{proof}
Applying the Taylor's formula, we have
\begin{align*}
\log(\rho_n^{2k_n})&=k_n\log(1+\rho_n^2-1)\\
&=k_n\Big(\rho_n^2-1-\frac{(\rho_n^2-1)^2}{2}+\ldots\Big).
\end{align*}
So $\rho_n^{2k_n}\to1$. Moreover, using the Taylor's formula again, we obtain, as $n\to\infty$
\begin{align*}
\rho_n^{2k_n}&=\exp\{k_n(\rho_n^2-1)+O(k_n(\rho_n^2-1)^2)\}\\
&=1+k_n(\rho_n^2-1)+O(k_n^2(\rho_n^2-1)^2).
\end{align*}
Thus the estimation (\ref{app0}) follows.
\end{proof}

\noindent\textbf{Proposition A.2.} {\it For model  (\ref{model}), assume that $y_0=0$ and $n(\rho_n-1)\to\gamma$ for some fixed $\gamma\neq0$ as $n\to\infty$, then we have
\begin{align*}
\frac{1}{n}\max_{1\leq i\leq n}|y_{i-1}|\rightarrow_p0,  \;\;\text{as}\;\; n\to\infty.
\end{align*}}

\begin{proof}
For any $\varepsilon>0$, by the Kolmogorov's maximal inequality,
\begin{align*}
&\mathbb{P}\Big(\frac{1}{n}\max_{1\leq i\leq n}|y_{i-1}|>\varepsilon\Big)=\mathbb{P}\Big(\max_{1\leq i\leq n}\big|\sum_{j=1}^{i-1}\rho_n^{i-1-j}\epsilon_j\big|>n\varepsilon\Big)\\
&\leq\frac{\sigma^2}{n^2\varepsilon^2}\cdot\sum_{j=1}^{n-1}\rho_n^{2(n-1-j)}
=\frac{\sigma^2(\rho_n^{2(n-1)}-1)}{n^2\varepsilon^2(\rho_n^2-1)}\to 0,
\end{align*}
as $n\to\infty$, which implies the desired result immediately. Here we used the facts, $n(\rho_n-1)\to\gamma$ and $\rho_n^n\to e^\gamma$ as $n\to\infty$.
\end{proof}

\begin{proof}[\textbf{Proof of Lemma \ref{l-s0}}]
(1). Note that there exists a sequence $m_n$ such that $m_n(1-\rho_n)\to\infty$ as $n\to\infty$, which implies that $\rho_n^{m_n}=o(1)$. So we can rewrite $\sqrt{1-\rho_n^2}y_0$ as follows
\begin{align*}
\sqrt{1-\rho_n^2}y_0&=\sqrt{1-\rho_n^2}\sum_{j=0}^{\infty}\rho_n^j\epsilon_{-j}\\
&=\sqrt{1-\rho_n^2}\sum_{j=0}^{m_n}\rho_n^j\epsilon_{-j}+\sqrt{1-\rho_n^2}\sum_{j=m_n+1}^{\infty}\rho_n^j\epsilon_{-j}\\
&=:\Gamma_{n,1}+\Gamma_{n,2}.
\end{align*}
Obviously, $\mathbb{E}(\Gamma_{n,2})=0$ and
\begin{align*}
\mathbb{E}(\Gamma_{n,2}^2)=\rho_n^{2(m_n+1)}\sigma^2=o(1),
\end{align*}
which immediately yield that $\Gamma_{n,2}\rightarrow_p0$ as $n\to\infty$. We apply the Lindeberg-Feller central limit theorem for the first term $\Gamma_{n,1}$. Denote $\zeta_{n,j}=\sqrt{1-\rho_n^2}\rho_n^{j}\epsilon_{-j}$, then $\mathbb{E}(\Gamma_{n,1})=0$ and
\[
\mathbb{E}(\Gamma_{n,1}^2)=(1-\rho_n^{2(m_n+1)})\sigma^2\to\sigma^2
\]
as $n\to\infty$. Hence, to complete the proof, we verify the Lindeberg condition,
\[
{\rm for~any~}\varepsilon>0,\qquad \sum_{j=0}^{m_n}\mathbb{E}\big(\zeta_{n,j}^2\mathbb{I}_{\{|\zeta_{n,j}|>\varepsilon\}}\big)=o(1).
\]
In fact, since $0<\rho_n<1$, by dominated convergence theorem, we have
\begin{align*}
\sum_{j=0}^{m_n}\mathbb{E}\big(\zeta_{n,j}^2\mathbb{I}_{\{|\zeta_{n,j}|>\varepsilon\}}\big)
&=(1-\rho_n^2)\sum_{j=0}^{m_n}\rho_n^{2j}\mathbb{E}\big(\epsilon_{-j}^2\mathbb{I}_{\{|\zeta_{n,j}|>\varepsilon\}}\big)\\
&\leq(1-\rho_n^2)\sum_{j=0}^{m_n}\rho_n^{2j}\mathbb{E}\big(\epsilon_{-j}^2\mathbb{I}_{\{\sqrt{1-\rho_n^2}|\epsilon_{-j}|>\varepsilon\}}\big)\\
&=(1-\rho_n^{2(m_n+1)})\mathbb{E}\big(\epsilon_1^2\mathbb{I}_{\{\sqrt{1-\rho_n^2}|\epsilon_{1}|>\varepsilon\}}\big)=o(1).
\end{align*}

(2). Because the Lindeberg's condition holds obviously in this case, we only need to estimate the asymptotic variance. Noticing that $\epsilon_t$ have zero median, we get $\mathbb{E}({\rm sign}(\epsilon_t))=0$ and
\[
\mathbb{E}\Big[\Big(\frac{1}{\sqrt{n}}\sum_{i=1}^n\rho_n^{i-1}{\rm sign}(\epsilon_i)\Big)^2\Big]=\frac{1-\rho_n^{2n}}{n(1-\rho_n^2)}\to1
\]
as $n\to\infty$, here we use the facts, $(1-\rho_n^n)\sim n(1-\rho_n)$ and $\rho_n^n\to1$ as $n\to\infty$.
\vskip5pt

(3). Since
\begin{align*}
\max_{1\leq i\leq n}|y_{i-1}|\leq\max\big\{|y_0|, \max_{2\leq i\leq n}|y_{i-1}|\big\}
\end{align*}
and
\begin{align*}
\max_{2\leq i\leq n}|y_{i-1}|\leq\max_{2\leq i\leq n}|\rho_n^{i-1}y_0|+\max_{2\leq i\leq n}\big|\sum_{j=1}^{i-1}\rho_n^{i-1-j}\epsilon_j\big|,
\end{align*}
part (1) of this lemma implies that
\begin{align*}
\sqrt{(1-\rho_n^2)/n}\max_{2\leq i\leq n}|y_{i-1}|&\leq\sqrt{(1-\rho_n^2)/n}\max_{2\leq i\leq n}\big|\sum_{j=1}^{i-1}\rho_n^{i-1-j}\epsilon_j\big|+o_p(1).
\end{align*}
By the Kolmogorov's maximal inequality and Proposition A.1, it follows that, for any $\varepsilon>0$,
\begin{align*}
&\mathbb{P}\Big(\sqrt{(1-\rho_n^2)/n}\max_{2\leq i\leq n}\big|\sum_{j=1}^{i-1}\rho_n^{i-1-j}\epsilon_j\big|>\varepsilon\Big)\\
&=\mathbb{P}\Big(\sqrt{(1-\rho_n^2)/n}\max_{2\leq i\leq n}\big|\sum_{j=1}^{i-1}\rho_n^{j-1}\epsilon_{i-j}\big|>\varepsilon\Big)\\
&\leq\frac{\sigma^2}{\varepsilon^2}\cdot\frac{1-\rho_n^{2n}}{n}=\frac{\sigma^2}{\varepsilon^2}\cdot\frac{1-\rho_n^{2n}}{n(1-\rho_n)}\cdot(1-\rho_n)=o(1),
\end{align*}
which complete the proof of part (3).
\end{proof}

\begin{proof}[\textbf{Proof of Lemma \ref{l-e0}}]
Noting that $\kappa_n(1-\rho_n)\to0$ and applying Proposition A.1, we can show that, $(1-\rho_n^{\kappa_n})\sim\kappa_n(1-\rho_n)$ and $\rho_n^{\kappa_n}\to1$ as $n\to\infty$. We omit the remainder of the argument since it is analogous to that in the proof of Lemma \ref{l-s0}.
\end{proof}

Proposition \ref{joint conv2} are required to achieve Theorem \ref{main result2} and \ref{main result3}. So, in the rest of this section, we mainly make the supplement for the proofs in Section \ref{s6}.  First of all, we need the following modified version of Theorem 1.4 in Chapter 7 of Ethier and Kurtz (1986).
\vskip5pt

\noindent\textbf{Theorem A.1} (Ethier and Kurtz, 1986) {\it Let $\Gamma(t)=(\gamma_{ij}(t))_{i,j=1}^d$ be a continuous, symmetric, $d\times d$ matrix-valued function, defined on $[0,1]$, satisfying $\Gamma(0)=0$ and for any $0\leq s<t\leq1$,
\[
\sum_{i,j=1}^d\big(\gamma_{ij}(t)-\gamma_{ij}(s)\big)\alpha_i\alpha_j\geq0,\qquad \forall\,\alpha=(\alpha_1,\ldots,\alpha_d)\in\mathbb{R}^d.
\]
Let $\{\xi_{k}^n,\mathcal{F}_{k}^n, n\geq1, 1\leq k\leq n\}$ be an $\mathbb{R}^d$-valued square-integrable
martingale difference array on a completed probability space $(\Omega,\mathcal{F},\mathbb{P})$.
Denote $S_n(t)=\sum_{k=1}^{[nt]}\xi_k^n$ and
\[
\Theta_n(t)=\big(\theta_n^{ij}(t)\big)_{i,j=1}^d=\sum_{k=1}^{[nt]}\mathbb{E}\big((\xi_k^n)^{\rm T}\xi_k^n|\mathcal{F}_{k-1}^n\big),
\]
 for $0\leq t\leq 1$. Suppose that,
\begin{align}\label{A1-1}
\mathbb{E}\Big(\sup_{0\leq t\leq 1}|S_n(t)-S_n(t-)|^2\Big)\to0, \quad as \; n\to\infty,
\end{align}
and for $i,j=1,2,\ldots,d$,
\begin{align}\label{A1-2}
\mathbb{E}\Big(\sup_{0\leq t\leq 1}|\theta_n^{ij}(t)-\theta_n^{ij}(t-)|\Big)\to0, \quad as \; n\to\infty.
\end{align}
In addition, if, for each $0\leq t\leq1$ and $i,j=1,2,\ldots,d$,
\begin{align}\label{A1-3}
\theta_n^{ij}(t)\rightarrow_p\gamma_{ij}(t), \quad as \; n\to\infty,
\end{align}
then the process $S_n$ converges weakly to a continuous process $G$ with independent Gaussian increments, mean vector zero and covariance matrix $\Gamma(t)$, on the Skorohod space $D_{\mathbb{R}^d}([0,1])$ of $\mathbb{R}^d$-valued cadlag paths on $[0,1]$.
}




\begin{proof}[\textbf{Proof of Lemma \ref{KLM1}}]
Recall that, $\mathcal{F}_i=\mathcal{F}_{n,i}=\sigma(\epsilon_1,\ldots,\epsilon_i)$ and $\mathcal{F}_0=\{\emptyset, \Omega\}$. For convenience, denote ${\bf X}_{n,i}=(K_{n,i},L_{n,i})$ and ${\bf X}_n(t)=\sum_{i=1}^{[nt]}{\bf X}_{n,i}$ for $0\leq t\leq 1$,
then it is clear that $\{{\bf X}_{n,i}, \mathcal{F}_i\}$ is an $\mathbb{R}^2$-valued square-integral martingale difference array.
Note that, under these circumstances,
\[
\Theta_n(t)=\sum_{i=1}^n\mathbb{E}\big({\bf X}_{n,i}^{\rm T}{\bf X}_{n,i}|\mathcal{F}_{i-1}\big)=\sum_{i=1}^n\mathbb{E}\big({\bf X}_{n,i}^{\rm T}{\bf X}_{n,i}\big).
\]
To apply Theorem A.1, we first verify the non-negative definiteness of the matrix-valued function $\Gamma(t)$ defined as in (\ref{covm}).
Clearly, for $0\leq s<t\leq 1$, the first principle minor of $\Gamma(t)-\Gamma(s)$ is positive for any $\gamma\neq0$; moreover, since $\sigma^2=1$, we know that $\mathbb{E}|\epsilon_1|\leq1$ and it follows that
\begin{align*}
\det(\Gamma(t)-\Gamma(s))=\frac{1}{4\gamma^2}\big(e^{\gamma(t-s)}-e^{-\gamma(t-s)}\big)^2-\big((t-s)\mathbb{E}|\epsilon_1|\big)^2\geq0.
\end{align*}
So the matrix-valued function $\Gamma(t)$ is non-negative definite.
\vskip5pt

Now, we check the condition (\ref{A1-1}), i.e.
\begin{align}\label{app-KLM1-0}
\mathbb{E}\Big(\max_{1\leq i\leq n}|{\bf X}_{n,i}|^2\Big)\to0,
\end{align}
as $n\to\infty$. Obviously, since $\rho_n^n\to e^{\gamma}$ as $n\to\infty$, we know that
\begin{align*}
\max_{1\leq i\leq n}|{\bf X}_{n,i}|^2\leq\frac{C_1}{n}\Big(1\vee\max_{1\leq i\leq n}\epsilon_i^2\Big),
\end{align*}
for some positive constant $C_1$ independent of $n$. For any $\varepsilon>0$, by the Chebyshev's inequality, we can get
\begin{align*}
\mathbb{P}(\max_{1\leq i\leq n}\epsilon_i^2>n\varepsilon)&\leq n\mathbb{P}(\epsilon_1^2>n\varepsilon)\\
&\leq\varepsilon^{-(2+\delta)/2}n^{-\delta/2}\cdot\mathbb{E}|\epsilon_1|^{2+\delta}\to0,
\end{align*}
as $n\to\infty$. So
\begin{align}\label{app-KLM1-1}
\max_{1\leq i\leq n}|{\bf X}_{n,i}|^2\rightarrow_p0.
\end{align}
On the other hand, notice that, for each $n$,
\[
\mathbb{E}\Big[\Big(\frac{\max_{1\leq i\leq n}\epsilon_i^2}{n}\Big)^{(2+\delta)/2}\Big]
\leq \frac{\mathbb{E}|\epsilon_1|^{2+\delta}}{n^{\delta/2}}.
\]
Hence the family of random variables, $\{\max_{1\leq i\leq n}|{\bf X}_{n,i}|^2, n\geq1\}$, is uniformly integrable, which, combined with (\ref{app-KLM1-1}), suggests that the condition (\ref{app-KLM1-0}) holds.
\vskip5pt

We next calculate the covariance matrix. Observe that
\begin{align*}
\sum_{i=1}^{[nt]}\mathbb{E}\big({\bf X}_{n,i}^{\rm T}{\bf X}_{n,i}\big)
=\left(
\begin{array}{cc}
\frac{1}{n}\sum_{i=1}^{[nt]}\rho_n^{2(i-1-n)} & t\rho_n^{-1}\mathbb{E}|\epsilon_1| \\
t\rho_n^{-1}\mathbb{E}|\epsilon_1| & \frac{1}{n}\sum_{i=1}^{[nt]}\rho_n^{2(n-i)} \\
\end{array}
\right).
\end{align*}
Therefore it is straightforward to show that, as $n\to\infty$,
\begin{align*}
\sup_{0\leq t\leq1}|\theta_n^{11}(t)-\theta_n^{11}(t-)|&=\frac{1}{n}\max_{1\leq i\leq n}\rho_n^{2(i-1-n)}\to0,\\
\sup_{0\leq t\leq1}|\theta_n^{22}(t)-\theta_n^{22}(t-)|&=\frac{1}{n}\max_{1\leq i\leq n}\rho_n^{2(n-i)}\to0,\\
\sup_{0\leq t\leq1}|\theta_n^{12}(t)-\theta_n^{12}(t-)|&=\sup_{0\leq t\leq1}|\theta_n^{21}(t)-\theta_n^{21}(t-)|=0,
\end{align*}
and
\[
\Theta_n(t)\rightarrow\Gamma(t)=\left(
                                  \begin{array}{cc}
                                    \frac{e^{-2\gamma}}{2\gamma}(e^{2\gamma t}-1) & t\mathbb{E}|\epsilon_1| \\
                                    t\mathbb{E}|\epsilon_1| & \frac{e^{2\gamma}}{2\gamma}(1-e^{-2\gamma t}) \\
                                  \end{array}
                                \right).
\]
Hence, by Theorem A.1, the process ${\bf X}_n$ converges weakly to a continuous process ${\bf{X}}=({\bf X}(t),0\leq t\leq1)$ which has independent Gaussian increments, zero mean vector and covariance matrix $\Gamma(t)$.

\vskip5pt

We now turn to the proof of the equation (\ref{KLM1-0}).
From the proof of Theorem 1.1 in Chapter 7 of Ethier and Kurtz (1986), we know that the limit process ${\bf X}$ can be represented by
\begin{align*}
{\bf X}^{\rm T}(t)=\int_0^t\Lambda^{1/2}(s)\,{\rm d}{\bf B}^{\rm T}(h(s)),\quad 0\leq t\leq1,
\end{align*}
where $\Lambda(\cdot)$ is a non-negative definite matrix-valued function, ${\bf B}=(B_1,B_2)$ is a $2$-dimensional
standard Brownian motion, and
\begin{align*}
h(t)=\gamma_{11}(t)+\gamma_{22}(t),\quad 0\leq t\leq1.
\end{align*}
In fact, the entries $\lambda_{ij}(\cdot)$ of $\Lambda(\cdot)$ are the solutions of the following equations,
\begin{align}\label{app-KLM1-2}
\left\{
  \begin{array}{llll}
    \gamma_{11}(t)=\int_0^t\lambda_{11}(s)\,{\rm d}h(s),   \\
    \gamma_{22}(t)=\int_0^t\lambda_{22}(s)\,{\rm d}h(s),   \\
    \gamma_{12}(t)=\int_0^t\lambda_{12}(s)\,{\rm d}h(s),   \\
    \gamma_{21}(t)=\int_0^t\lambda_{21}(s)\,{\rm d}h(s).
  \end{array}
\right.
\end{align}
Solving (\ref{app-KLM1-2}), we can get
\begin{align*}
\lambda_{11}(t)&=\frac{e^{-2\gamma(1-t)}}{e^{2\gamma(1-t)}+e^{-2\gamma(1-t)}}=\frac{1-\tanh(2\gamma(1-t))}{2},\\
\lambda_{22}(t)&=\frac{e^{2\gamma(1-t)}}{e^{2\gamma(1-t)}+e^{-2\gamma(1-t)}}=\frac{1+\tanh(2\gamma(1-t))}{2},\\
\lambda_{12}(t)=&\lambda_{21}(t)=\frac{\mathbb{E}|\epsilon_1|}{e^{2\gamma(1-t)}+e^{-2\gamma(1-t)}}=\frac{\mathbb{E}|\epsilon_1|}{2\cosh(2\gamma(1-t))}.
\end{align*}

Hence we complete the proof of Lemma \ref{KLM1}.
\end{proof}

\noindent\textbf{Theorem A.2} (Rootz\'{e}n, 1980) {\it Let $(B(t), 0\leq t\leq1)$ be a standard Brownian motion with respect to the filtration $(\mathcal{F}_t, 0\leq t\leq1)$. Suppose $(\phi_n(t), 0\leq t\leq 1)_{n\geq1}$ is a sequence of random  functions which is adapted to the filtration $(\mathcal{F}_t, 0\leq t\leq 1)$. If
\begin{align*}
\sup_{0\leq t\leq 1}\left|\int_0^t\phi_n(s)\,{\rm d}s\right|\rightarrow_p0, \quad {as} \; n\to\infty,
\end{align*}
and
\begin{align*}
\int_0^1\phi_n^2(t)\,{\rm d}t\rightarrow_p\tau,\quad {as} \; n\to\infty,
\end{align*}
for some random variable $\tau$ such that $\tau>0$ a.s., then
\begin{align*}
\frac{\int_0^1\phi_n(t)\,{\rm d}B(t)}{\sqrt{\int_0^1\phi_n^2(t)\,{\rm d}t}}\rightarrow_d\mathscr{N}(0,1), \quad {as} \; n\to\infty.
\end{align*}
}
\vskip5pt

\noindent\textbf{Theorem A.3} (Graversen and Peskir, 2000) {\it Let $(V(t),  t\geq0)$ be the Ornstein-Uhlenbeck process solving
 \[
 {\rm d}V(t)=-\beta V(t){\rm d}t+{\rm d}B(t)
 \]
 with $V(0)=0$, where $\beta>0$ and $(B(t),t\geq0)$ is a standard Brownian motion, then there exist universal positive constants $c_0,c_1$ such that, for all stopping times $\tau$ of $(V(t),t\geq0)$,
 \begin{align*}
 \frac{c_0}{\sqrt{\beta}}\mathbb{E}\sqrt{\log(1+\beta\tau)}\leq\mathbb{E}\left(\max_{0\leq t\leq \tau}|V(t)|\right)\leq\frac{c_1}{\sqrt{\beta}}\mathbb{E}\sqrt{\log(1+\beta\tau)}.
 \end{align*}
}

\begin{proof}[\textbf{Proof of Lemma \ref{result3-lem1}}]
Since
\[
\sup_{0\leq t\leq 1}\left|\int_0^t\varphi_\gamma(s)\,{\rm d}s\right|\leq\int_0^1|\varphi_\gamma(t)|\,{\rm d}t,
\]
it is enough to show $\int_0^1|\varphi_\gamma(t)|\,{\rm d}t\rightarrow_p0$ as $\gamma\to\infty$. Notice that, from Lemma \ref{KLM1},
\[
L(t)\sim\mathscr{N}\Big(0,\frac{e^{2\gamma}(1-e^{-2\gamma t})}{2\gamma}\Big),
\]
so
\begin{align*}
\mathbb{E}\int_0^1|\varphi_\gamma(t)|\,{\rm d}t&=2\gamma e^{-\gamma}\int_0^1e^{-\gamma(1-t)}\mathbb{E}|L(t)|\,{\rm d}t\\
&=\sqrt{4\gamma/\pi}e^{-\gamma}\int_0^1\sqrt{e^{2\gamma t}-1}\,{\rm d}t\\
&\leq\sqrt{4/(\pi\gamma)}(1-e^{-\gamma})\to0, \quad as \; \gamma\to\infty.
\end{align*}
By the Markov's inequality, this establishes the equation (\ref{result3-lem1-1}).
\vskip5pt

We next prove equation (\ref{result3-lem1-2}). Firstly, using integration by parts, we can get
\begin{align*}
\int_0^1\varphi_\gamma^2(t)\,{\rm d}t&=4\gamma^2e^{-4\gamma}\int_0^1e^{2\gamma t}L^2(t)\,{\rm d}t\\
&=2\gamma e^{-4\gamma}\int_0^1 L^2(t)\,{\rm d}e^{2\gamma t}\\
&=2\gamma e^{-4\gamma}\left(L^2(1)e^{2\gamma}-\int_0^1e^{2\gamma t}\,{\rm d}L^2(t)\right).
\end{align*}
Furthermore, since
\[
\widetilde{\lambda}_{12}^2(t)h'(t)+\widetilde{\lambda}_{22}^2(t)h'(t)=e^{2\gamma(1-t)}, \quad 0\leq t\leq1,
\]
the It\^{o}'s formula, together with equation (\ref{main result3-1}), yields that
\[
{\rm d}L^2(t)=2L(t)\,{\rm d}L(t)+e^{2\gamma(1-t)}\,{\rm d}t.
\]
Therefore
\begin{align*}
\int_0^1\varphi_\gamma^2(t)\,{\rm d}t=2\gamma e^{-2\gamma}L^2(1)-2\gamma e^{-2\gamma}-4\gamma e^{-4\gamma}\int_0^1e^{2\gamma t}L(t)\,{\rm d}L(t).
\end{align*}
Noticing that
\[
L(1)\sim\mathscr{N}\Big(0,\frac{e^{2\gamma}-1}{2\gamma}\Big),
\]
we have
\begin{align*}
2\gamma e^{-2\gamma}L^2(1)=_d\big(\sqrt{1-e^{-2\gamma}}\mathscr{N}(0,1)\big)^2\rightarrow_p\mathscr{N}^2(0,1),
\end{align*}
as $\gamma\to\infty$. In addition,
\begin{align*}
\mathbb{E}\left(4\gamma e^{-4\gamma}\int_0^1e^{2\gamma t}L(t)\,{\rm d}L(t)\right)^2&=16\gamma^2e^{-6\gamma}\int_0^1e^{2\gamma t}\mathbb{E}L^2(t)\,{\rm d}t\nonumber\\
&=8\gamma e^{-4\gamma}\int_0^1(e^{2\gamma t}-1)\,{\rm d}t\leq4e^{-2\gamma}\to0,
\end{align*}
as $\gamma\to\infty$. It follows from the Chebyshev's inequality that
\[
4\gamma e^{-4\gamma}\int_0^1e^{2\gamma t}L(t)\,{\rm d}L(t)\rightarrow_p0, \quad as \; \gamma\to\infty.
\]
Hence we achieve the equation (\ref{result3-lem1-2}). The proof is complete.
\end{proof}

\begin{proof}[\textbf{ Proof of Lemma \ref{result3-lem2}}]
We first show that the equation (\ref{result3-lem2-1}) holds. Let
\begin{align*}
{\rm d}{B}^*(t)=\frac{\widetilde{\lambda}_{12} (t)\sqrt{h'(t)}{\rm d}\widetilde{B}_1(t)+\widetilde{\lambda}_{22} (t)\sqrt{h'(t)}{\rm d}\widetilde{B}_2(t)}{\sqrt{\widetilde{\lambda}_{12}^2(t)h'(t)+\widetilde{\lambda}_{22}^2(t)h'(t)}},
\end{align*}
then, by the L\'{e}vy characterization of Brownian motion, $B^*$ is a $1$-dimensional standard Browinan motion. Since
\[
\widetilde{\lambda}_{12}^2(t)h'(t)+\widetilde{\lambda}_{22}^2(t)h'(t)=e^{2\gamma(1-t)}, \quad 0\leq t\leq1,
\]
the equation (\ref{main result3-1}) implies that
$$
{\rm d}L(t)=e^{\gamma(1-t)}{\rm d}B^*(t).
$$
Hence, by the It\^{o}'s formula, we have
\begin{align*}
\int_0^t\psi_\gamma(s)\,{\rm d}s&=\sqrt{-2\gamma}e^{-\gamma}\int_0^te^{\gamma s}L(s){\rm d}s\\
&=\frac{\sqrt{-2\gamma}e^{-\gamma}}{\gamma}\left(e^{\gamma t}L(t)-\int_0^t e^{\gamma s}\,{\rm d}L(s)\right)\\
&=\frac{\sqrt{-2\gamma}}{\gamma}\left(\int_0^te^{\gamma(t-s)}\,{\rm d}B^*(t)-B^*(t)\right).
\end{align*}
Therefore
\begin{align*}
\sup_{0\leq t\leq1}\left|\int_0^t\psi_\gamma(s)\,{\rm d}s\right|\leq\sqrt{\frac{2}{-\gamma}}
\left(\sup_{0\leq t\leq 1}\left|\int_0^te^{\gamma(t-s)}\,{\rm d}B^*(t)\right|
+\sup_{0\leq t\leq 1}|B^*(t)|\right).
\end{align*}
Obviously,
\[
\sqrt{\frac{2}{-\gamma}}\left(\sup_{0\leq t\leq 1}|B^*(t)|\right)\to_p0,
\]
as $\gamma\to-\infty$ and Theorem A.3 shows that
\begin{align*}
\mathbb{E}\left(\sqrt{\frac{2}{-\gamma}}
\sup_{0\leq t\leq 1}\Big|\int_0^te^{\gamma(t-s)}\,{\rm d}B^*(t)\Big|\right)
\leq\frac{c_1\sqrt{2\log(1-\gamma)}}{-\gamma}\to0,
\end{align*}
as $\gamma\to-\infty$. By the Markov's inequality again, we achieve the equation (\ref{result3-lem2-1}).
\vskip5pt

Finally, we verify the equation (\ref{result3-lem2-2}). Notice that
\[
{\rm d}L^2(t)=2L(t)\,{\rm d}L(t)+e^{2\gamma(1-t)}\,{\rm d}t.
\]
It follows that
\begin{align*}
\int_0^1\psi_\gamma^2(t)\,{\rm d}t=-L^2(1)+1+2
e^{-2\gamma}\int_0^1e^{2\gamma t}L(t)\,{\rm d}L(t).
\end{align*}
Since
\[
L(t)\sim\mathscr{N}\Big(0,\frac{e^{2\gamma}(1-e^{-2\gamma t})}{2\gamma}\Big),
\]
we have
$$
L^2(1)\to_p0, \quad as\; \gamma\to-\infty.
$$
Moreover,
\begin{align*}
\mathbb{E}
\left(2e^{-2\gamma}\int_0^1e^{2\gamma t}L(t)\,{\rm d}L(t)\right)^2
&=4e^{-4\gamma}\int_0^1e^{4\gamma t}e^{2\gamma(1-t)}\mathbb{E}L^2(t)\,{\rm d}t\\
&=\frac{2}{\gamma}\int_0^1(e^{2\gamma t}-1)\,{\rm d}t\to0,
\end{align*}
as $\gamma\to-\infty$. It concludes from the Chebyshev's inequality that
\[
2e^{-2\gamma}\int_0^1e^{2\gamma t}L(t)\,{\rm d}L(t)\to_p0.
\]
Hence we obtain the equation (\ref{result3-lem2-2}).  The proof is complete.
\end{proof}

\section*{Acknowledgment}

The authors would like to express their sincere gratitude to
the anonymous referees and AE for helpful comments which surely lead to an improved presentation
of this paper. The authors are very grateful to Huarui He, Hui Jiang, Feng Li, Yu Miao, Shaochen Wang and Qingshan Yang for the helpful discussions. Hailin Sang's work was partially supported by the Simons Foundation grant 586789, USA. Guangyu Yang's work was partially supported by the Foundation of Young Scholar of the Educational Department of Henan Province grant 2019GGJS012, China.

\vskip20pt

\noindent {\it Data Availability Statement:} The data that support the findings of this study are available from the corresponding author upon reasonable request.

\end{document}